\documentclass[11pt]{article}

\usepackage{pstricks}
\usepackage[all]{xypic}
\usepackage[ansinew]{inputenc}
\usepackage[psamsfonts]{amssymb}
\def\hH{{H\!H}} 
\def\hC{{\bf C}} 
\def\cT{{T\! C}} 
\def\bB{I\!\!B} 

\title{Rational String Topology}
\author{Yves F\'elix, Jean-Claude Thomas and Micheline Vigu\'e-Poirrier }

\begin{document}
\maketitle

    \begin{abstract} We use the   computational power of rational homotopy theory to
  provide   an explicit cochain  model for the loop product and the string bracket
of a 1-connected closed manifold $M$. We   prove that the loop
homology of $M$ is isomorphic to the Hochschild cohomology of the
commutative graded algebra $A_{PL}(M)$ with coefficients in
itself.  Some explicit computations of the loop product and the
string bracket
 are given.
\end{abstract}

    \vspace{5mm}\noindent {\bf AMS Classification} : 55P35, 54N45,55N33,
  17A65,
    81T30, 17B55

    \vspace{2mm}\noindent {\bf Key words} : String homology, rational homotopy,
Hochschild cohomology, free loop space, loop homology.

 \section{Introduction and preliminaries}

Let $M$ be a 1-connected closed  oriented $m$-manifold and
$$
LM=M^{S^1}:=\{
\mbox{unbased continuous map }   S^1 \to M\}
$$
 be the associated free loop space.  The
{\sl loop homology}  of $M$ is the ordinary  homology of $LM$,
with a shift of degrees, ${\mathbb H}_*(LM)=H_{*+m}(LM)$ together
with a commutative associative product, called {\em loop product}:
$$
{\mathbb H}_p(LM) \otimes {\mathbb H}_q(LM) \to {\mathbb
H}_{p+q}(LM)\,, \quad a\otimes b \mapsto a \bullet b\,.
$$
This product was discovered by M. Chas and D. Sullivan, \cite{CS}
who defined it in terms of ``transversal geometric chains'' on
$LM$ by mixing the intersection product on $M$, $$ H_p(M) \otimes
H_q(M) \to H_{p+q- m}(M) $$ with the composition of loops $$
LM\times _MLM =\{\mbox{composable loop}\}  \to LM \,.$$  In
\cite{CJ}, R. Cohen and J.D.S. Jones used the Pontryagin-Thom
structure to show that the loop product is realized at the
homotopy level by  Thom spectra. By using intersection product in
the setting of Hilbert manifolds, \cite{Ch},  D. Chataur, provides
an other description of the loop product.

In this paper we start with an elementary description of the loop
product   in terms of the diagonal class of $M$.   Then, we prove

\vspace{3mm}\noindent {\bf  Theorem A.} {\sl  When coefficients
are taken in a field of characteristic zero, there exists an
explicit model of the loop product in terms of the minimal model
of Sullivan of $M$.}

\vspace{3mm} This model is described in details in the text.
Theorem A implies that the rational loop product does not depend
on the geometry of the manifold but only on its rational homotopy
type.

 When coefficients are in a field, there is also an
isomorphism of graded vector spaces, \cite {CJ}, \cite{FTV} from
the singular homology  of the space $LM$ to  the Hochschild
cohomology of the singular cochain complex of $M$:  $$
  {\mathbb H}_\ast (LM) \to \hH ^{\ast} (C^\ast  M, C^\ast  M) \,.
$$

In \cite{CJ}, Cohen and Jones assert that this isomorphism can be
modified into an isomorphism of commutative graded algebras.  More
recently  S. A. Merkulov, \cite{M}, gives a proof based on the
theory of iterated integrals when the coefficients are the real
numbers.

We construct here a direct isomorphism of algebras when the
coefficients are the rational numbers. Our construction uses
differential graded Lie algebras.

 Recall that the Lie model of the space $M$, denoted $L$, is a differential graded Lie algebra
 defined by the property that
 the cochain algebra $C^*L$ is quasi-isomorphic to the Sullivan model of
 $M$. In particular,
if $C_\ast(L)$ denotes the chain complex of the differential
graded Lie algebra $L$ then $H\ast(C_\ast(L)=H_\ast(M;{\mathbb Q}
$.

   The key point of our work
 is a chain morphism connecting the chain and the cochain complex of $L$ with coefficients
  in $(UL)_a^\vee = \mbox{Hom}(UL,\mathbb Q)$ with the adjoint representation:
 $$C^*(L;(UL)^\vee_a) \stackrel{cap}{\longrightarrow} C_{m-*}(L;(UL)^\vee_a)\,.$$
The morphism $cap$ is the cap product with a cycle in $C_m(L)$ representing
 the fundamental class in   homology.
 Thanks to the diagonal  on  $UL$, the cochain complex
 $C^* (L;(UL)^\vee_a)$ is a commutative differential graded  algebra
  quasi-isomorphic to the Sullivan minimal model of $LM$. On the other hand, thanks
 to the
 multiplication  in $UL$, the chain complex $C_*(L;(UL)^\vee_a)$ is a differential graded
 coalgebra,
 and the cohomology of the dual differential graded algebra   is isomorphic to the Hochschild cohomology
 $H\!H^*(UL,UL)$. More precisely we   prove:

 \vspace{3mm}\noindent {\bf Theorem B.}  {\sl Let $L = {\cal L}_M$ be the minimal Lie model
 of $M$.
 Then the dual map $\mbox{Hom}(cap, \mathbb Q)$ induces in homology an isomorphism  of algebras
 $$  {\mathbb H}_*(LM) \stackrel{\cong}{\leftarrow} H^*(C^*(L;(UL)_a)) \,,$$
 where $(UL)_a$ is equipped with the adjoint representation and
 $H_*(LM)$ with the loop multiplication.
 }

 \vspace{3mm}\noindent
Now using homological algebra arguments, we prove that
$H^*(C^*(L;(UL)_a)) $ is isomorphic to the Hochschild cohomology $
H\!H^*(UL,UL) $. Then using the natural isomorphism
$$H\!H^*(UL,UL)\cong H\!H^*(C^*L,C^*L)\,,$$ (\cite{FMT}),
 and the existence, for the minimal Lie model $L$ of $M$, of quasi-isomorphisms of
 algebras,
  $C^*L \stackrel{\cong}{\leftarrow} A\stackrel{\cong}{\rightarrow} C^*M\,,$  we recover

 \vspace{3mm}\noindent {\bf Theorem C.} {\sl Let $M$ be a simply connected closed manifold,
 then we have a natural isomorphism of algebras
 $\Phi : {\mathbb H}_*(LM) \stackrel{\cong}{\to} H\!H^*(C^*M,C^*M)$.}

 \vspace{3mm} In \cite{CS}, Chas and Sullivan construct also a
 morphism of algebras $$I : {\mathbb H}_*(LM)\to H_{*}(\Omega M)$$  defined
 as follows.
 Let $x_0$ be the base point  of $M$ and $[x_0]\in H_0(LM)$ be the homology class of the constant
loop, then for any $a\in H_*(LM)$, we define $I(a) = a \bullet
[x_0] $. The map $I$ looks like an intersection with the fiber
$\Omega M$.
  We prove

 \vspace{3mm}\noindent {\bf Theorem D.}  {\sl There exists an isomorphism of
 algebras $\Psi$
 making commutative the diagram
 $$
 \begin{array}{ccc}
   {\mathbb H}_*(LM) &\stackrel{ \Phi }{\longrightarrow}& H\!H^*(C^*M,C^*M)\\
   {\scriptstyle I}\downarrow && \downarrow {\scriptstyle H\!H^*(C^*M,\varepsilon)}\\
   H_*(\Omega M) & \stackrel{\Psi}{\longrightarrow} & H\!H^*(C^*M;\mathbb Q)
   \end{array}$$
   where $\varepsilon : C^*M\to \mathbb Q$ denotes the usual augmentation.}

\vspace{2mm} Let us recall that an element $x \in \pi_q(M)$ is
called a {\it
 Gottlieb element} (\cite{FHT}-p.377), if the map
 $x \vee id_M : S^q \vee M \to M$ extends to the product $S^q \times M$.
 The Gottlieb
elements
 generate a  subgroup $G_\ast (M)$ of $\pi_\ast (\Omega M)$ via
 the isomorphism $\pi_\ast (\Omega M) \cong \pi_{\ast +1}(M)$. Finally, we denote
by
 cat $M$ the Lusternik-Schnirelmann category of $M$, normalized so
 that cat\,$S^n$=1.  Ii follows from (\cite{FTV}-Theorem 2)  that :

 \vspace{3mm}

{\it  a)  The kernel of $I$ is a nilpotent ideal of nilpotency
index less than or equal to cat $M $.

 b) (Im $  I$) $ \cap\, (\pi_\ast (\Omega M)\otimes {\mathbb Q}
 ) = G_\ast ( M)\otimes {\mathbb Q}$.

 c) $\displaystyle\sum_{i=0}^n$ dim (Im $   I)_i  \leq C n^k $,  some constant $C>0$ and $k
 \leq$ cat $M $. }

\vspace{4mm} Our model of the loop product can also be used to
compute the string bracket. Denote by $LM\times_{S^1}ES^1$ the
equivariant free loop space. The long exact homology sequence
associated to the sphere bundle $S^1 \to LM\times
ES^1\stackrel{p}{\to} LM\times_{S^1}ES^1$ has the form $$\to
H_n(LM)\stackrel{H_n(p)}{\to} H_n(LM\times_{S^1}ES^1)
\stackrel{Cap}{\to}
 H_{n-2}(LM\times_{S^1}ES^1)
\stackrel{M}{\to} H_{n-1}(LM)\to \cdots$$ where $Cap$ is the cap
product with the characteristic class of the sphere bundle. Using
this sequence, Chas and Sullivan give  to
$H_*(LM\times_{S^1}ES^1)$ the structure of
 a graded Lie algebra
of degree $(2-m)$
$$[a,b]= (-1)^{\vert a\vert} H_*(p)( M(a)\bullet M(b))\,,\quad \mbox{(\cite{CS})}\,.$$

Few things are known about this bracket. For surfaces of genus larger than zero, they recover
formulas proved in the context of symplectic geometry.

Using computations made in \cite{BV} we prove that for any
1-connected manifold $M$ such that the cohomology algebra
$H^*(M,\mathbb{Q})$ is monogenic, the string bracket is trivial.
So for spheres, complex projective spaces..., string bracket  is
trivial. We also produce  a model allowing explicit computations
of the string bracket   for any 1-connected compact manifold.

 \vspace{3mm} The text is divided in two parts. In the first one, after recalling some materials,
  we give   a description
 of the loop  product and the string bracket in terms of the Sullivan minimal model. This step
 makes the bridge between topology (geometry) and algebra. In the
 second part, we construct an explicit algebra isomorphism $\mathbb
 H_*(LM) \cong H\!H^*(C^*M;C^*M)$, and we prove
Theorems B, C  and D.

 \vspace{4mm}
Part 1 - Algebraic models

\vspace{2mm}

2. Preliminaries on differential homological algebra

3. A 'Sullivan model' description of the loop product

4. The string bracket

\vspace{3mm} Part 2 - The isomorphism $\mathbb
 H_*(LM) \cong H\!H^*(C^*M;C^*M)$

\vspace{2mm}

5. Hochschild cohomology of a Lie algebra

6. Lie models and the cap product $C^*(L;M) \to C_{m-*}(L;M)$

7. Proof of Theorem B.

8. The intersection morphism

\section{Preliminaries on differential homological algebra}

\vspace{3mm} All the graded vector spaces, algebras, coalgebras and Lie algebras $V$ are
defined over $\mathbb Q$ and are supposed of finite type, i.e. dim $V_n <\infty$ for all $n$.

\subsection{Graded vector spaces}

If $V =\{V_i\}_{i \in \mathbb Z}  $ is a
(lower)
 graded $\mathbb Q$-vector space  (when we need upper graded vector space  we put
 $V_i= V^{-i}$ as usual) then  $V^{\vee}$ denotes the graded  dual vector space,
    $$V^{\vee} = \mbox{Hom}(V,\mathbb Q)\,,$$
$sV$ denotes the suspension of $V$, $$(sV)_n =V_{n-1}\,, \,\, (sV)^n =V^{n+1}\,,$$
and  $TV $ denotes the tensor algebra on $V$,
while we denote by $T\!C(V)$ the free
 supplemented coalgebra generated by $V$.

 Since we work with
    graded   objets, we will make a
    special attention to signs.  Recall that if
  $P=\{P_i\}$ and  $N= \{N_i\}$ are differential
graded vector spaces  then
\begin{enumerate}\item[$\bullet$]  $P\otimes N$ is a differential graded vector space :
$$(P\otimes N)_r =
\oplus _{p+q=r} P_p
  \otimes N_q\,, \quad d_{P\otimes N}= d_P \otimes id_N + id_P
  \otimes d_N,.$$
\item[$\bullet$]  $\mbox{Hom}(P,N)$ is a differential graded vector space  : $$Hom_n(P,N)=
\prod_{k-l=n} \mbox{Hom}(P_l,N_k)\,,
 \quad D_{Hom(P,N)}f  = d_N\circ f - (-1) ^{|f|} f\circ d_P\,.$$
 \end{enumerate}

  If $C$ is a differential graded coalgebra with diagonal
$\Delta$
 and $A$ is a differential graded algebra with product $\mu$ then the
 cup product,
$f\cup g= \mu\circ (f\otimes g) \circ \Delta$, gives  to the differential
graded vector space
$\mbox{Hom}(C,A)$ a structure of differential graded  algebra.

\subsection{Two-sided bar construction}

  Let  $(A, d)$ be a differential graded supplemented
algebra,
    $ A = \mathbb Q \oplus \overline A$, and $(P, d)$ and $(N,d)$  be
     differential graded $ A$-modules, respectively   right and  left $A$-module.
    The {\it two-sided bar construction}, $\bB(P; A; N)$  is defined as follows:
    $$
    \bB _k (P; A; N)  =  P \otimes T^k
    (s\overline A ) \otimes N\,.$$
     A generic element   is written $
    p[a_1|a_2|...|a_k]n$  with  degree   $| p|+ | n|+ \sum_{i=
    1}^k (|s a_i|)\,.$ The differential $d= d_0+d_1$ is defined by
    $$
    d_0  : \bB_k  (P;A;N) \to \bB _k  (P;A;N)\,, \quad  d_1  : \bB _k
    (P;A;N) \to \bB_{k-1}  (P;A;N)\,,
    $$
    with
    $$
    d_0  p[\, ] n = d (p) [\, ] n + (-1)^{| p|} m [\, ]d(n) \,, \qquad
    d_1   p[\,]n=0
    $$
    and if $k>0$, and
    $\epsilon _i = | p| + \sum _{j<i} (|s a_j|)$:
    $$
    \renewcommand{\arraystretch}{1.5}
    \begin{array}{rll}
    d_0  ( p[a_1|a_2|...|a_k]n)& = d( p) [a_1|a_2|...|a_k]n  -
   \displaystyle\sum
    _{i=1}^k (-1)^{\epsilon _i}  p[a_1|a_2|...|d(
    a_i)|...|a_k]n\\ &+ (-1) ^{\epsilon _{k+1}}  p[a_1|a_2|...|a_k]d(n)
    \\[2mm]
    d_1  p[a_1|a_2|...|a_k]n &= (-1) ^{| p|}  pa_1[a_2|...|a_k]n +
    \displaystyle\sum _{i=2}^k (-1) ^{\epsilon _i}
  p[a_1|a_2|...|a_{i-1}a_i|... |
   a_k]n \\
    &- (-1)^{\epsilon _{k}}  p[a_1|a_2|...|a_{k-1}]
    a_k n
    \end{array}
     \renewcommand{\arraystretch}{1}
    $$

    \vspace{3mm}

The complex
  $
 \bB A ={\bB} (\mathbb Q;A;\mathbb Q)
 $ is a differential graded coalgebra whose comultiplication is defined by
 $$\Delta [a_1\vert \cdots \vert a_r] =
 \sum_{i=0}^r [a_1\vert \cdots \vert a_i]\otimes [a_{i+1}\vert \cdots   \vert a_r]\,.$$

\subsection{ Hochschild cohomology}

  \vspace{3mm}   Let $ A$ be a differential
 graded algebra and $N$ a differential graded $A$-bimodule.
 The {\it Hochschild cochain complex} of $A$ with coefficients in   $M$ is    the   cochain  complex
    $$
    {\hC}^\ast  (A;N) =  (\mbox{Hom}   (\cT(sA), N), D_0+D_1)
    = \mbox{Hom}_{A^{e}}(\bB(A;A;A),N)\,,
    $$
    with $A^e = A \otimes A^{opp}$.
The differential $D_0+D_1$ comes from the natural differential on
$\mbox{Hom}_{A^{e}}(\bB(A;A;A),N)$:   $D_0(f)([\,]) =
d_N(f([\,]))$, $D_1(f)([\,])=0$, and for $k\geq 1$ and
  $f \in  \mbox{Hom}( T\!C(sA), N) $,  then
   $$D_0(f)([a_1|a_2|...|a_k])
  = d_{N}(f([a_1|a_2|...|a_k])) + (-1)^{|f|} \sum _{i=1} ^k
  (-1)^{\overline \epsilon_i} f([a_1|...|d_Aa_i|...|a_k])$$  and
  $$\renewcommand{\arraystretch}{1.6}
\begin{array}{ll}
D_1(f)([a_1|a_2|...|a_k])= &
   - (-1)^{|f|} \left( (-1)^{|a_1|\, |f|}  a_1 f([a_2|...|a_k])\right)\\ &
(-1)^{|f|} \left(
  \sum _{i=2} ^k
  (-1)^{\overline \epsilon_i} f([a_1|...|a_{i-1}a_i|...|a_k])\right)\\
 & (-1)^{|f|} \left(
  + (-1) ^{\overline\epsilon _k} f([a_1|a_2|...|a_{k-1}])a_k
  \right)\,,
\end{array}
\renewcommand{\arraystretch}{1}
$$
 where
  $\overline \epsilon _i = |sa_1|+|sa_2|+...+|sa_{i-1}|$.

\vspace{2mm}\noindent
 The cohomology of $\hC ^\ast(A;N)$ is called the  {\it Hochschild cohomology of $A$ with coefficients in $N$} is
and is denoted by  $
 \hH^*(A;N) \,.
 $

The natural cup product on $\mbox{Hom}(\cT(s\bar A),A)$ defined by
$$f\cup g = \mu_{A} \circ (f\otimes g) \circ \Delta_{T\!C(s\bar A)}\,,$$ makes $\hC^*(A;A)$
 a   differential graded algebra, and
 $\hH^* (A;A) $ into a commutative graded algebra (\cite{GS}).

\subsection{The chain coalgebra and the cochain algebra on a differential graded Lie algebra}

The {\it chain coalgebra} on a differential graded Lie algebra $(L,d_L)$
is the graded differential coalgebra
$$C_*L =(\land sL,d_0+d_1)\,.$$
The differential and the comultiplication are defined by
$$d_0(sx_1\land \cdots \land sx_k) = -\sum_{i=1}^k (-1)^{\sum_{j<i}\vert sx_j\vert}\,\,
sx_1\land \cdots \land sd_Lx_i\land \cdots \land sx_k\,,$$
$$d_1(sx_1\land \cdots \land sx_k) = \sum_{1\leq i<j\leq k}
(-1)^{e_{ij}}\,\, s[x_i,x_j]\land \cdots \widehat{sx_i}  \cdots
\widehat{sx_j} \cdots \land sx_k\,,$$
$$\Delta (sx_1\land \cdots \land sx_k) = \sum_{j=0}^k \,\sum_{\sigma\in Sh(j )}
\varepsilon_\sigma (sx_{\sigma (1)} \land \cdots \land sx_{\sigma (j)})
 \otimes (sx_{\sigma (j+1)} \land \cdots \land sx_{\sigma (k)})\,,$$
 where $\varepsilon_\sigma$ is the usual sign associated to the graded permutation
 $\sigma$, $Sh(j )$ denotes the set of $(j,k-j)$-shuffles, and
 $$e_{ij}=   \vert sx_i\vert (1 + \sum_{k<i}\vert sx_k\vert)
 \, + \, \vert sx_j\vert (\sum_{k<j, k\neq i} \vert sx_k\vert ) \,.$$

\vspace{3mm} For each left $L$-module $N$ and right $L$-module $P$ we can define the complex
$$C_*(P;L;N) =(P \otimes C_*L \otimes N, d_0+d_1)\,,$$
with
$$d_0(p\otimes c \otimes n) = dp\otimes c \otimes n + (-1)^{\vert p\vert }
p\otimes d_0 (c)\otimes n  + (-1)^{\vert p\vert + \vert c\vert}
p\otimes c \otimes dn\,,$$
$$
\renewcommand{\arraystretch}{1.5}
\begin{array}{ll}
d_1(p\otimes sx_1 \land &\cdots \land sx_k\otimes n) =   \\
&
(-1)^{\vert p\vert} p \otimes d_1(sx_1\land \cdots\land sx_k)\otimes dn\\
& + \sum_{i=1}^k (-1)^{\vert p\vert + (\sum_j\vert sx_j\vert) +
\vert x_i\vert + \vert sx_i\vert (\sum_{j>i} \vert sx_j\vert)}
p\otimes sx_1 \land\cdots \widehat{sx_i} \cdots \land sx_k\otimes x_i\cdot n\\
& + \sum_{i=1}^k (-1)^{\vert p\vert + \vert sx_i\vert (\sum_{j<i} \vert sx_j\vert)}
p\cdot x_i \otimes sx_1 \land \cdots \widehat{sx_i} \cdots \land sx_k\otimes n\,.
\end{array}
\renewcommand{\arraystretch}{1}
$$

The chain complex of $L$ with coefficients in a left-module $N$,
$$C_*(L;N) := C_*(\mathbb Q;L;N)\,$$
 is a left
$C_*L$-comodule. In a similar way   the chain complex $C_*(P;L) :=
C_*(P;L;\mathbb Q)$ is a $C_*L$-right comodule. Here $\mathbb Q$
is equipped with the trivial action.

When $N = UL$ with the action induced by left multiplication then $C_*(L, UL)$ is a
left $C_*L$-comodule and a right $UL$-module and both structures are compatible.
Moreover the augmentation
$$  C_*(L;UL)= C_*L \otimes UL
\stackrel{\varepsilon\otimes \varepsilon}{\longrightarrow}  \mathbb Q$$
is a quasi-isomorphism.

The {\it cochain complex} of $L$ with coefficients in a right
$L$-module $R$ is defined by
 $$C^*(L;R) = \mbox{Hom}_{UL} (C_*(L;UL), R)\,.$$

\vspace{2mm} The homology and the cohomology of $L$ with coefficients in a module $N$ are defined by
$H_*(L;N) = H_*(C_*(L;N))$ and $H^*(L;N) = H^*(C^*(L;N))$.

\vspace{2mm} The {\it cochain  algebra} on a differential graded
Lie algebra $(L,d_L)$ is the graded differential  algebra $C^*L =
C^*(L;\mathbb Q) = \mbox{Hom}(C_*L,\mathbb Q)$.

\subsection{Semifree resolutions}

\vspace{3mm}  Let $A$ be a differential graded algebra. A module
$P$ is called {\sl a semifree module} if $P$ is equipped with a
filtration $P = \cup_{n\geq 0} P(n)$, satisfying $P(0) = 0$, $P(n)
\subset P(n+1)$ and such that $P(n)/P(n-1)$ is free on a basis of
cycles (\cite{FHT}).

For any  $A$-module $N$,  there exists a semifree module $P$ and a
quasi-isomorphism
 $\varphi : P \to N$. The module $P$ is called a semifree resolution of $N$.

For instance the complex
$\bB(A;A;A)\buildrel{\simeq}\over\rightarrow A$ is a semifree
resolution of $A$ as an $A^e$-module~\cite[4.3(ii)]{FHT}.

The complex $C_*(UL;L;UL) = UL \otimes C_*L\otimes UL$ is a semifree resolution of $UL$
as an $UL$-bimodule. This follows from the quasi-isomorphism of $UL$-bimodules
 $\theta : C_*(UL;L;UL) \to UL$  defined by
 $\theta (a\otimes 1 \otimes b) = ab$ and $\theta
(a\otimes c\otimes b) = 0$ for $a, b \in UL$ and $c
\in \land^+(sL)$.

\section{A 'Sullivan model' description of the loop product}

\subsection{A definition  of the loop product at the homological level}

Let $M$ be a closed riemannian $m$-dimensional manifold and
$\Delta M \hookrightarrow M\times M $ be the diagonal embeding of
$M$ in   $M\times M$. We denote by $\nu$ the normal bundle of
$\Delta M$ in $M\times M$ which is identified with the tangent
bundle $\tau M$ of $M$. The exponential map of the normal bundle
$\nu$,   defined only on a neighborhood of the zero section of
$\nu$, induces a diffeomorphism from a $m$-dimensional disk
bundle,   $\tau_DM$,  onto a tubular neighborhood, denoted $N$, of
$\Delta M$.   The associated $(m-1)$-sphere bundle, $\tau _SM$ is
diffeomorphic to the boundary $\partial N$ of $N$.

 Let $ p:
L(M\times M)= LM\times LM \to M\times M\, \quad p(\omega, \omega')
= (\omega (0), \omega'(0))$ be   the product of the two loop
fibrations. We consider the commutative diagram

\xymatrix{
 LM\times_MLM =p^{-1}(\Delta M) \ar[r]^{\mbox{} \hspace{1cm}j'}\ar[d]^q  & p^{-1}(N) \ar[r]^{i'}\ar[d]^{p'}
   & LM\times LM \ar[d]^p &
  p^{-1}( M\times M-\Delta M)\ar[d]^{p''}\ar[l]_{k'}\\
 M=\Delta M \ar[r]^{j}  & N \ar[r]^{i}  & M\times M
&   M\times M-\Delta M\,, \ar[l]_{k} }

\vspace{2mm}\noindent where $i,j,k,i',j'$ and $k'$ are the natural
inclusions.

 One   observes that in
the  above diagram  the vertical maps $p',q$ and $p''$ are
pullback   fibrations and that
 $p^{-1}(N)$ can be identified
 with
 the pullback of the fiber bundle  $\tau_DM$ along $q$: $ p^{-1}(N) \cong q^*\tau_DM\,.$

The loop multiplication on $H_*(LM)$ with coefficients in $\mathbb
Q$ is then defined by the composition of maps

\xymatrix@R=5mm { H_*(LM) \otimes H_*(LM)\ar[d]_{\cong} &&\\
  H_*(LM\times LM)\ar[r] &
H_*(LM\times LM, p^{-1}(M\times M - \Delta
M))\ar[d]^{\cong}_{E}&\\
 & H_*(p^{-1}(N), p^{-1}(\partial N))\ar[d]_{\cong}
&\\
 & H_*(q^*\tau_DM, q^*\tau_SM)\ar[d]_{Th_*} &\\
& H_{*-m}(LM\times_MLM)\ar[r]^{H_{*-m}(\gamma)}& H_{*-m}(LM)}

\vspace{2mm}\noindent Here $\gamma : LM\times_MLM \to LM$ denotes
the composition of loops, $E$ is the classical excision
isomorphism and $Th_*$ is the Thom isomorphism recalled below.

We recall that, if $p : X \to M$ is a $m$-dimensional vector
bundle with associated disk bundle,  $X_D$, and sphere bundle
$X_S$, there exists a class $\alpha \in H^m(X_D,X_S)$, called the
Thom class of the disk bundle, such that the cap product with
$\alpha$ followed by the projection on $M$ induces an isomorphism
$Th_* : H_*(X_D,X_S) \to H_{*-m}(M)$. Here we consider the normal
bundle $\nu$ defined above and its $m$-dimensional disk bundle
$\tau_DM$, that defines the Thom class $\alpha \in
H^m(\tau_DM,\tau_SM)$. The isomorphism $Th_*:
H_*(q^*\tau_DM,q^*\tau_SM)\to H_{*-m}(LM\times_MLM)$ is the cap
product with $q^*\alpha$ followed by the projection.

 \subsection{The Sullivan minimal model of $M$}

 The Sullivan minimal model of a  simply connected closed $m$-dimensional manifold $M$,
${\cal M}_{M}$, is a commutative differential graded algebra
defined over the rational numbers that represents the rational
homotopy type of $M$  (\cite{S},\cite{FHT}).
 As an algebra ${\cal M}_{M}$ is a free graded commutative algebra,
 $${\cal M}_{M} = \land V\,,$$
 with $V^n \cong \mbox{Hom}(\pi_n(M),\mathbb Q)$. The differential $d$ satisfies
 the
  minimality condition $d(V) \subset \land^{\geq 2}(V)$.

  The minimal model is unique  up to isomorphism. It depends only on the rational homotopy
  type of $M$ and contains all the rational homotopy informations of $M$. In particular
   there exists a
  differential graded algebra $B$ and quasi-isomorphisms
  $${\cal M}_{M} \stackrel{\simeq}{\leftarrow} B \stackrel{\simeq}{\to} C^*(M;\mathbb Q)\,.$$
There exists also a quasi-isomorphism of differential graded
algebras $\psi : {\cal M}_{M} \otimes  \mathbb R \to \Omega_{DR}
M$, where $\Omega_{DR} M$ denotes the algebra of de Rham forms on
$M$.

  More generally a  Sullivan (or free) model for $M$ is a commutative differential graded algebra
  $(\land W,d)$, that is free as a commutative graded algebra,
   for which there exists a quasi-isomorphism $\varphi : {\cal M}_{M} \to (\land W,d)$,
   with $W = W^{\geq 2}$ a finite type graded vector space.
  It is sometimes more convenient to work with free models having special properties than
   to work with minimal models, as we will see further.

\subsection{Examples of models}

  Let $(\land V,d)$ be a free model of $M$. A free model for the free loop space $LM$
  is   given by
   $(\land V \otimes \land  sV, D)$, where $(sV)^n = V^{n+1}$, $D(sv)=-s(dv)$, where $s$
   has been extended to
   $\land V \otimes \land sV$ as a derivation satisfying $s(v) = sv$ and $s(sv) = 0$
   (\cite{VS}).

A model of $LM\times_MLM$ is then given by $(\land V \otimes \land
sV\otimes \land s'V,D)$ where $(s'V)^n = V^{n+1}$ and $D(s'v) =
-s'(dv)$. In other words, $(\land V\otimes \land s'V,D)$ is a copy
of $(\land V\otimes \land sV,D)$ and $$(\land V\otimes \land sV
\otimes \land s'V,D)  = (\land V\otimes \land sV,D) \otimes_{\land
V}  (\land V\otimes \land s'V,D)\,.$$

A model for the diagonal map $\Delta : M \to M\times M$ is given
by the multiplication $\mu : (\land V,d)\otimes (\land V',d) \to
(\land V,d)$ where $(\land V',d)=(\land V,d)$. The Sullivan
relative model of $\mu$,   is a   commutative diagram in which
$\varphi$ is a quasi-isomorphism and $i$ the canonical inclusion

\vspace{2mm}\hspace{2cm} \xymatrix{ (\land V,d)\otimes (\land
V',d) \ar[r]^{\mu}\ar[rd]^{i} & (\land V,d) \\  & (\land V \otimes
\land V' \otimes \land \bar V,D)\,.\ar[u]_{\varphi}
 }

\vspace{3mm}\noindent   Here $\bar V^n = V^{n+1}$ and $D(\bar v) -
(v-v') \in \bar V \oplus \land^{\geq 2} (V \oplus V' \oplus \bar
V) $.

\subsection{A Sullivan model for the loop product}

In the next section, we will construct an explicit model for the
composition of paths $LM\times_MLM \to LM$ that has the form
$$c : (\land V \otimes \land sV,D) \to (\land V \otimes \land
sV\otimes \land s'V,D)\,.$$

Denote now by $ u_M $ a cycle in $\mbox{Hom}(\land V,\mathbb Q)$
representing the fundamental class in homology:  $[u_M]\in
H_m(M;\mathbb Q)= H_m(\mbox{Hom}(\land V,\mathbb Q))$. We choose
an homogeneous basis $\alpha_i$ of $H^*(M)$ and its Poincar\'e
dual basis $\alpha_i^{\#}$:
$$\langle \alpha_i \cup\alpha_j^{\#} ; [u_M]\rangle =
\delta_{ij}\,.$$ We choose then  cocycles $a_i$ and $a_i^{\#}$ in
$\land V$ whose cohomology classes are respectively $\alpha_i$ and
$\alpha_i^{\#}$. This gives  the cocycle   $$T = \sum_i
(-1)^{\vert\, a_i\vert} a_i \otimes a_i^{\#}\in \land V \otimes
\land V\,.$$ By Theorem 1 below, the multiplication by $T$,
$$\mu_T : (\land V,d)\otimes (\land V',d) \to (\land V,d)\otimes
(\land V',d)$$ extends into a morphism of $(\land V\otimes \land
V')$ differential modules $$\mu_T : (\land V \otimes \land V'
\otimes \land \bar V,D) \to (\land V,d)\otimes (\land V',d)\,.$$
We have thus defined all the terms involved in the diagram

\hspace{1cm}\xymatrix{ (\land V \otimes \land sV,D)\ar[d]^c  &\\
  (\land V \otimes (\land sV \otimes \land s'V),D) &
 (\land V \otimes \land  V' \otimes \land \bar V  \otimes (\land s V\otimes \land s' V),D')
 \ar[l]_{ \varphi\otimes 1}\ar[d]^{ \mu_T\otimes 1}  \\   & (\land V
\otimes \land sV, D)\otimes (\land V' \otimes \land s'V, D)
 }

\vspace{3mm}\noindent   Theorem A of the Introduction can be now
expressed in   more explicit terms

\vspace{3mm}\noindent {\bf Theorem A.}  {\sl The induced map in
cohomology $$H^*(\mu_T\otimes 1)\circ H^*(\varphi\otimes 1)^{-1}
\circ H^*(c) : H^{*-m}(LM) \to H^*(LM) \otimes H^*(LM)$$ is the
coproduct on $H^*(LM)$ dual to the loop product on $H_*(LM)$. }

\vspace{3mm}\noindent {\bf Proof.} We have to compute the map
induced in cohomology by the composite (1)
  $$
\begin{array}{ccc}
 C^{*-m}(LM\times_MLM)  &&\\  {\scriptstyle \cup
\beta}{\downarrow}&&\\
  C^*(q^*\tau_DM, q^*\tau_SM)
 \stackrel{\cong}{\to} & C^*(p^{-1}(N),
p^{-1}(\partial N))&  \stackrel{E}{\to}  C^*(LM\times LM, p^{-1}(M
\times M - \Delta M))\\ &&\downarrow\\ && C^*(LM\times LM)
    \end{array}
$$ where $\beta$ is the Thom class of the disk bundle $p^{-1}(
N)\to LM\times_MLM$.  Since this fiber bundle is the inverse image
of the disk bundle $   N \to M$, $\beta = q^*(\beta')$, where
$\beta'$ is the Thom class of $\tau_DM$. We observe that the
diagram (1) is obtained by pullbacking $p$ along the maps arising
in the diagram $$
\begin{array}{ccc}
 &C^{*-m}(M) &\\
 (2) &   {\scriptstyle \cup \beta'}{\downarrow}&\\
 &  C^*( \tau_DM,  \tau_SM)&
 \stackrel{\cong}{\to}   C^*(N,
  \partial N)     \stackrel{
E}{\to}  C^*( M\times  M,
 (M \times M- \Delta M))\to C^*( M\times  M) \,.
    \end{array}
$$ A   model for the composite (1), as a $(\land V\otimes \land
sV)\otimes (\land V'\otimes\land s'V)$-module is thus obtained by
making the tensor product of a model of (2) with  $(\land V\otimes
\land sV)\otimes (\land V'\otimes\land s'V)$ over $\land V\otimes
\land V'$. Now it follows from (\cite{Br}, 11.2) that the map
induced in cohomology by the composite (2) is the cohomology
intersection, i.e. the multiplication by $[T]$. Therefore by
Theorem 1 below, a model of
 (1) is given by    $\mu_T $. This
ends the proof of Theorem A. \hfill $\square$

\vspace{3mm} \noindent{\bf Theorem 1.}  {\sl There exists  up to
homotopy one and only one   morphism of differential graded
$(\land V\otimes \land V')$-module $$ f : (\land V\otimes \land
V'\otimes \land\bar V,D) \to (\land V,d) \otimes (\land V',d)$$
whose restriction to $\land V\otimes \land V'$ gives in homology
the multiplication by $T$.}

\vspace{3mm}\noindent {\bf Proof.}  The commutation with the
differentials gives to the graded vector space $$\mbox{Hom}_{\land
V\otimes \land V'}(\land V\otimes \land V'\otimes \land \bar V,
\land V\otimes \land  V')$$ a differential whose homology classes
in degree $d$ are the homotopy classes of module morphisms of
degree $d$. By definition of $\mbox{Ext}$, we have,
$$H^*(\mbox{Hom}_{\land V\otimes \land V'}(\land V\otimes \land
V'\otimes \land \bar V, \land V\otimes \land V') =
\mbox{Ext}_{\land V\otimes \land V'} ^*(\land V,\land V \otimes
\land V')$$ where $\land V$ is viewed as a $\land V\otimes \land
V'$-module by the multiplication map.

We use the Moore spectral sequence and consider only the total
degree $$\mbox{Ext}^*_{H\otimes H}(H,H\otimes H) \Rightarrow
\mbox{Ext}_{\land V\otimes \land V'} ^*(\land V,\land V\otimes
\land V')\,,$$ with $H = H(\land V,d)$.

Since $H\otimes H$ is a Poincaré duality algebra,
$\mbox{Ext}^q_{H\otimes H}(\mathbb Q,H\otimes H) = 0$ for $q \neq
2m$ and dim\,$\mbox{Ext}^{2m}_{H\otimes H}(\mathbb Q,H\otimes H)=
1$. Therefore by induction on the dimension, we have that
$$\mbox{Ext}_{H\otimes H}(E,H\otimes H)=\mbox{Ext}^{\geq
2m-d}_{H\otimes H}(E,H\otimes H)$$ when $E$ is finite dimensional
and $E = E^{\leq d}$. Since $M$ is simply connected, $H^{m-1}=0$,
and   the long exact  sequence associated to the short exact
sequence $ 0 \to H^m\to H \to H/H^m\to 0$ gives
$\mbox{Ext}_{H\otimes H}^m(H,H\otimes H) \cong\mathbb Q\alpha$,
for some $\alpha$, and $\mbox{Ext}_{H\otimes H}^p(H,H\otimes H) =
0 $ for $p=m-1$ and $p=m+1$. Therefore the element $\alpha$ will
be an $\infty$-cycle and never a boundary. In particular
$$\mbox{Ext}_{\land V\otimes \land V'} ^m(\land V,\land V\otimes
\land V')\cong \mathbb Q \,.$$

Since the multiplication by $T$,  $H \to H\otimes H$ is a morphism
of $H\otimes H$-bimodules, this multiplication is the representant
of a generator of $\mbox{Ext}_{H\otimes H}^m(H, H\otimes H)$.

This proves the existence and the unicity of a map $f$ that
extends the multiplication by $T$ on $\land V\otimes \land V'$.
\hfill$\square$

\vspace{3mm} Remark that the existence of   $\mu_T$ and the
construction of the model of the composite can also be deduced
from a more general construction given by Lambrechts and Stanley
in \cite{LS}.

\vspace{3mm}\subsection{Construction of a model for $LM\times_MLM \to LM$}

The free loop space $LM$ is the  pullback of the diagram $$
\begin{array}{ccc}
LM & \to &M^{[0,1]}\\
\downarrow && \mbox{}\hspace{5mm}\downarrow {\scriptstyle (p_0,p_1)} \\
M & \stackrel{\Delta}{\to} & M\times M\,,
\end{array}
$$
where $p_i(\omega) =\omega (i)$.
In the same way, the space $LM\times_MLM$ is the pullback of the diagram
$$
\begin{array}{ccc}
LM\times_MLM & \to &M^{[0,1]}\times_MM^{[0,1]}\\
\downarrow && \downarrow {\scriptstyle q} \\
M & \stackrel{\Delta}{\to} & M\times M\times M\,,
\end{array}
$$ where $q(\omega, \omega') = (\omega(0), \omega(1)=\omega'(0),
\omega'(1))$. The composition of paths $LM\times_MLM \to LM$ is
the pullback of $id_M$ and the   path composition $\nu$ over
$\rho$:

\begin{pspicture}(-2,0)(10,5)
 \psset{xunit=6mm}
    \psset{yunit=6mm }
\psline{->}(4.3,0)(9,0)
\psline{->}(4.4,5)(9,5)\psline{->}(2.5,7)(5,7)\psline{->}(1.3,2)(5.7,2)
\psline{->}(4,4.7)(4,0.4)
\psline{->}(1,6.5)(1,2.4)\psline{->}(7,6.7)(7,2.4)\psline{->}(1.3,1.8)(3.7,0.2)
\psline{->}(1.6,6.6)(3.7,5.2)\psline{->}(7.6,6.6)(9.7,5.2)\psline{->}(7.3,1.8)(9.7,0.2)
\psline{->}(10,4.6)(10,0.4) \rput(7,7){$\scriptstyle
M^{[0,1]}\times_MM^{[0,1]}$} \rput(1,7){$\scriptstyle
LM\times_MLM$} \rput(1,2){$M$} \rput(7,2){$\scriptstyle M\times
M\times M$}
 \rput(7.2,4){$\scriptstyle q$}
 \rput(5,2.3){$\scriptstyle \Delta$}
\rput(4,5){$\scriptstyle LM$} \rput(4,0){$\scriptstyle M$}
\rput(10,0){$\scriptstyle M\times M$} \rput(10,5){$\scriptstyle
M^{[0,1]}$} \rput(7,0.3){$\scriptstyle \Delta$}
\rput(2.8,1.4){$\scriptstyle id_M$} \rput(10.9,2.5){$\scriptstyle
(p_0,p_1)$} \rput(8.7,1.2){$\scriptstyle
\rho$}\rput(8.6,6.2){$\scriptstyle \nu$}
\end{pspicture}

\vspace{3mm}\noindent Here $\rho (\alpha,\beta,\gamma) =
(\alpha,\gamma)$.

 Denote now by $\sigma : M \to M^{[0,1]}$
the map sending a point   to the constant path at that point. Then
  $\sigma$ is an homotopy
 equivalence making  commutative the diagram

\vspace{3mm} \mbox{}\hspace{3cm}
\begin{pspicture}(0,0)(11,4.5)
 \psset{xunit=6mm}
    \psset{yunit=6mm }
\psline{->}(6.7,4)(3.8,4)\psline{->}(3.5,6)(1.1,6)
\psline[linestyle=dashed]{->}(3.7,5.7)(0.3,2.3)
\psline{->}(3,3.7)(3,0.4)
\psline{->}(0,5.7)(0,2.4)\psline{->}(6.8,3.8)(3.4,0.4)
\psline{->}(0.6,5.6)(2.7,4.2)\psline{->}(4.3,5.8)(6.7,4.2)\rput(5.5,5.6){$\scriptstyle
id_M$}\psline{->}(0.3,1.8)(2.4,0.4) \rput(5.1,4.2){$\scriptstyle
\sigma$} \rput(-1,6){$\scriptstyle M^{[0,1]}\times_MM^{[0,1]}$}
\rput(0,2){$\scriptstyle M\times M\times
M$}\rput(4,6){$\scriptstyle M$}
 \rput(0.2,4){$\scriptstyle q$}\rput(5,2.5){$\scriptstyle\Delta$}
 \rput(1.5,3.9){$\scriptstyle\Delta$}\rput(2.5,6.4){$\scriptstyle(\sigma,\sigma)$}
\rput(3,4){$\scriptstyle M^{[0,1]}$} \rput(3,0){$\scriptstyle
M\times M$}\rput(7,4){$\scriptstyle M$}
  \rput(1.5,1.3){$\scriptstyle\rho$}
\rput(3.8,2.4){$\scriptstyle(p_0,p_1)$}
\rput(1.5,5.3){$\scriptstyle\nu$}
\end{pspicture}

\vspace{4mm}\noindent  Therefore a model of $\nu$ is   a model of
$id_M$. A model of the loop composition $LM\times_MLM \to LM$  is
then obtained by the tensor product over a model of $\rho$ of a
model of $id_M$ and a relative model of $id_M$.  Denote by

\vspace{2mm}\hspace{2cm} \xymatrix{ (\land V \otimes \land V'
\otimes \land sV,D) \ar[r]^{\mbox{}\hspace{14mm}\varphi}  & (\land
V,d)\\ (\land V\otimes \land V',D) \ar[u]^{i} \ar[ur]_{ \mu}& }

\vspace{3mm} \noindent a model of the multiplication $\mu$ as
defined in 3.3: $(\land V',d)$ is a copy of $(\land V,d)$, $(sV)^n
= V^{n+1}$, $D(sv) - v+v' \in \land^{\geq 2}
 (V \oplus V' \oplus sV)\,\oplus sV$, $\mu(v) = \mu(v') = v$ and $\varphi(sv)=0$.

 We
first construct a relative model $c'$ of $\nu$.

\vspace{3mm} \mbox{}\hspace{3cm}
\begin{pspicture}(0,0)(11,4)
 \psset{xunit=6mm}
    \psset{yunit=6mm }
\psline{->}(5,4)(6.7,4)\psline{->}(1.1,6)(3.5,6)
\psline[linestyle=dashed]{-}(0.3,2.3)(1.7,3.7)
\psline[linestyle=dashed]{->}(2.6,4.6)(3.7,5.7)
\psline{->}(3,0.4)(3,3.5)
\psline{->}(0,2.4)(0,5.7)\psline{->}(3.2,0.2)(6.8,3.8)
\psline{->}(2.7,4.2)(0.6,5.6)\psline{->}(6.7,4.2)(4.3,5.8)\rput(5.5,5.6){$\scriptstyle
id_{\land V} $}\psline{->}(2.7,0.2)(0.6,1.6)
\rput(5.1,4.3){$\scriptstyle\varphi$} \rput(0,6){$\scriptstyle Z$}
\rput(0,2){$\scriptstyle\land V\otimes \land V'\otimes \land
V''$}\rput(4,6){$\scriptstyle\land V$}
 \rput(0.2,4){$\scriptstyle i$}\rput(5,2.3){$\scriptstyle\mu$}\rput(2.8,5.1){$\scriptstyle\mu$}
 \rput(2.5,6.4){$\scriptstyle\varphi\otimes_{\land V}\varphi$}
\rput(3,4){$\scriptstyle\land V\otimes \land V'\otimes \land sV$}
\rput(3,0){$\scriptstyle\land V\otimes \land V'
$}\rput(7,4){$\scriptstyle\land V$}
  \rput(1.5,1.3){$\scriptstyle\gamma$}
\rput(3.2,2.4){$\scriptstyle i$}  \rput(1.5,5.4){$\scriptstyle
c'$}
\end{pspicture}

\vspace{4mm}\noindent Here $Z = (\land (V \oplus V'\oplus
V''\oplus  sV \oplus s'V),D)$, $\gamma(a\otimes b') = a\otimes
b''$, $(\land V'\otimes \land V'' \otimes \land  s'V,D)$ is a copy
of $ (\land V \otimes \land V'\otimes  \land sV ,D)$ and $$(\land
(V \oplus V'\oplus V''\oplus  sV \oplus s'V),D) = (\land (V \oplus
V' \oplus  sV )
 ,D)\otimes_{\land V'} ( \land (V'\oplus V''\oplus s'V),D)\,.$$

A model $c$ for the path composition is then  obtained
  by the tensor product of $c'$ with $(\land V,d)$ over $\land V\otimes \land V'\otimes \land V''$.
$$
\begin{array}{ccc}
(\land V\otimes \land V' \otimes \land sV,D) & \stackrel{c'}{\to}
& (\land (V \oplus V'\oplus V''\oplus  sV \oplus s'V),D)\\
{\scriptstyle \pi}\downarrow && {\scriptstyle \pi}\downarrow\\
(\land V \otimes \land sV, D ) & \stackrel{c}{\to} & (\land V
\otimes \land sV \otimes \land s'V,D)
\end{array}
$$
Here $\pi(v)=\pi(v') = \pi(v'')= v$, $\pi(sv) = sv$ and $\pi(s'v)=s'v$.

\subsection{The example of formal spaces}

\vspace{3mm}  Let $M$ be a formal space, i.e. a space $M$ whose
minimal model $(\land V,d)$ is quasi-isomorphic to $(H^*(M),0)$.
  Examples of formal spaces
are given by simply connected compact K\"ahler manifolds
(\cite{DGMS}) and homogeneous spaces of rank zero
 (quotient of
compact connected Lie groups by closed subgroups of the same
rank).

Denote by $(\land V\otimes \land \bar V,D)$  the minimal model of
the free loop space $LM$. When $M$ is a formal space the coproduct
dual to the loop product is the map in cohomology induced by the
composition $$(H^*(M)\otimes \land \bar V,D) \stackrel{c}{\to}
(H^*(M)\otimes \land ( \bar V \oplus \bar V'),D) \stackrel{
\mu_T\otimes1 }{\longrightarrow} (H^*(M)\otimes \land \bar V,D)
\otimes (H^*(M)\otimes \land \bar V,D) \,.$$ with $(H^*(M)\otimes
\land \bar V, D)= H^*(M)\otimes_{\land V}(\land V\otimes \land
\bar V,D)$, and where $c$ is a model for the path composition
$LM\times_MLM \to LM$.

\vspace{2mm} \noindent {\bf The particular case $M = \mathbb
CP^{n}$.}  This is a formal space whose  minimal model is given by
$(\land (x,y),d)$,
 $d(y) = x^{n+1}$, $\vert x\vert = 2$, $\vert y\vert = 2n+1$.

 The model of the free loop space is
 $$(\land (x,\bar x, y, \bar y), D)\,,\quad D(\bar x) = 0\,,
  D(\bar y) = -(n+1)x^{n}\bar x\,.$$  Since we have a
  quasi-isomorphism $(\land (x,\bar x, y, \bar y), D)\to (\land
  (x,\bar x, \bar y)/(x^{n+1}),d)$,
  the cohomology of the free loop space is
 $$H^*(\land (x,\bar x,   \bar y)/(y^{n+1}), D) \cong \mathbb Q\cdot 1 \oplus  \left(\land^+(x,\bar x)
 /(x^{n+1}, x^n\bar x) \otimes \land  \bar y \right)\,.$$

A model $c$ of the path composition $  LM\times_MLM \to LM$  is
given by
 $$c (\bar x) = \bar x + \bar x'\,,\hspace{1cm}
 c (\bar y) = \bar y+\bar y' - \frac{n(n+1)}{2} x^{n-1}\bar x\bar x'\,.$$
 The dual   of the loop product is quite easy to handle
  and is induced by the   map
 $$\theta : (H^*(M)\otimes \land  (\bar x,  \bar y),D) \to
 (H^*(M)\otimes \land  (\bar x,  \bar y),D) \otimes (H^*(M)\otimes \land
 (\bar x, \bar y),D),$$
 $$
 \renewcommand{\arraystretch}{2}
 \left\{
 \begin{array}{ll}
 \theta (\alpha \otimes  \,  \bar y^{[s]} ) =& \displaystyle\sum_{p=0}^{n}
 \sum_{j=0}^s
\,\, \alpha x^p    \bar y^{[j]} \otimes x^{n-p} \bar y^{[s-j]}\\&
- \displaystyle\frac{n(n+1)}{2}  \sum_{p=0}^{n}
 \sum_{j=0}^{s-1}\,\, \alpha x^{n-1+p}\bar x    \bar y^{[j]}
\otimes x^{n-p}\bar x \bar y^{[s-j]} \,,\\ \theta (\alpha\otimes
\bar x \otimes \bar y^{[s]} ) =& (1\otimes \bar x + \bar x\otimes
1) \cdot \left(\displaystyle\sum_{p=0}^{n}
 \sum_{j=0}^s
\,\, \alpha x^p    \bar y^{[j]} \otimes x^{n-p}
y^{[s-j]} \right)\,, \end{array}\right.
\renewcommand{\arraystretch}{1}$$ with $\alpha\in H^*(M)$ and
$\bar y^{[s]} =  \bar y^s / s!$.

We consider the   basis for the free loop space homology formed by
the elements $$1\,,a_{p,q}\,, b_{r,s}\,, \quad   \,
 p=1, \ldots , n,\,\,\, q\geq 0,\,\,\, s\geq 0,\,\,\, r=0, \ldots ,n-1 \,,$$
 with $\vert a_{p,q}\vert = 2p+ 2qn $,
$\vert b_{r,s}\vert = 2r+1+ 2sn $,
 corresponding by
the isomorphism $H^*(LM) = \mbox{Hom}(H_*(LM),\mathbb Q)$ to the
classes $1$, $x^p\bar y^{[q]}$, and $x^p\bar x\bar y^{[q]}$.
$$\langle a_{p,q}, x^r\bar y^{[s]}\rangle =\left\{
\begin{array}{ll} 1 & \mbox{if}\,\, (p,q)=(r,s)\\ 0 &
\mbox{otherwise}
\end{array}
\right.\,, \quad
 \langle b_{p,q}, x^r\bar x \bar y^{[s]}\rangle =\left\{ \begin{array}{ll}
1 & \mbox{if}\,\, (p,q)=(r,s)\\ 0 & \mbox{otherwise}
\end{array}
\right.$$
 From the description of the map $\theta$, we deduce that the loop product is described by the
formulas $$ a_{p,q}\bullet a_{r,s} = a_{p+r-n, q+s}\,,
\hspace{5mm} a_{p,q}\bullet b_{r,s} = b_{p+r-n,q+s}\,,
\hspace{5mm} (a_{n-1,0})^n = 1\,, \hspace{5mm} 1\bullet a_{n,1} =
0\, .$$ This shows that $$\mathbb H_*(L(\mathbb CP^n);\mathbb Q)
\cong \land(a,b, t)/(a^{n+1}, a^nb, a^nt)\,,$$ with $\vert a \vert
= -2$, $\vert b\vert = -1$ and $\vert t \vert = 2n$,\,
$a=a_{n-1,0}$,\, $b=b_{n-1,0}$,\, $t=a_{n,1}$. Note that this
computation can also be found  in \cite{CJY} as an application of
a spectral sequence defined by Cohen, Jones and Yan.

\section{The string bracket}

Following Chas and Sullivan \cite {CS}, we call string homology
the equivariant homology of the free loop space $${\mathcal H}_* =
H_{*+m}^{S^1}(LM)= H_{*+m}(LM\times_{S^1} ES^1;\mathbb Q)\,.$$ The
circle fibration $ S^1 \to LM\times ES^1 \stackrel{p}{\to}
LM\times_{S^1}ES^1$ leads to the
 exact Gysin sequence
$$\cdots \rightarrow  \mathbb H_n \stackrel{
H(p)}{\longrightarrow} {\mathcal H}_n
\stackrel{c}{\longrightarrow} {\mathcal H}_{n-2} \stackrel{\mathbb
{M}}{\longrightarrow} \mathbb H_{n-1}\rightarrow \cdots $$ where
$c$ is the cap product with the characteristic class of the circle
bundle.

\vspace{3mm}\noindent {\bf Definition.}  The string bracket on
${\mathcal H}_*$ is the binary
 operation$[-,-]: {\cal H}_*\otimes {\cal H}_* \rightarrow  {\cal H}_{*-2} $ defined by
$$[a,b]= (-1)^{\vert a\vert} H_*(p) (\mathbb {M}(a) \bullet
\mathbb {M}(b))\,,$$ where the product $\bullet$  is the loop
product on $\mathbb H_*(LM)$. Chas and Sullivan prove that
$({\mathcal H}_*, [-,-])$ is a graded Lie algebra of degree $2$
(\cite{CS}, Theorem 6.1).

\vspace{3mm} Using results on the rational cohomology of
$LM\times_{S^1} ES^1$, \cite {BV}, and the splitting of the long
exact Gysin sequence in cohomology for formal spaces, \cite {V},
it is easy to prove:

\vspace{3mm} \noindent {\bf Theorem 2.}  {\sl Let $M$ be a
1-connected closed manifold such that $H^*(M,\mathbb Q)$ is a
truncated algebra in one generator, then the string bracket on the
string homology $\mathcal H_*$ is trivial.}\\

\vspace{3mm}\noindent {\bf Proof.}

 From \cite {BV}, we have the following facts:\\
If $H^*(M) = \Lambda u$ with $\vert u\vert = 2p+1$ then $\mathcal
H_{2i} = 0$ for all $i$.\\ If $H^*(M) = \Lambda u/(u^{n + 1})$
with $\vert u\vert = 2p$, then $\dim_{\mathbb Q}\mathcal H_{2i} =
1$ for all $i$. Furthermore, the space $M$ is formal and it is
shown in \cite {V} that the map $c$
 is an isomorphism.\\
In the two cases above, we get the nullity of the maps: $$ \mathbb
E = H_*(p): \mathbb H_{2i}\longrightarrow \mathcal H_{2i} $$ $$
\mathbb {M}: \mathcal H_{2i}\longrightarrow \mathbb H_{2i+1} $$
Let $a\in \mathcal H_{2i-1} $ and $b\in \mathcal H_{2b-1}$, for
some $(i,j)\in \mathbb Z^2$, then $\mathbb M(a)\bullet \mathbb
M(b) \in \mathbb H_{2(i + j +1)}$, so we have $$ [a,b] = -\mathbb
E(\mathbb M(a)\bullet \mathbb M(b)) = 0. $$ \hfill$\square$

 \vspace{3mm}\noindent{\bf Description of the Gysin sequence.}  When $(\land V,d)$ denotes the
minimal model of $M$, the Gysin sequence in cohomology can be
easily computed using models (\cite{VS}, \cite{BV}). First of all
a model of $LM$ is $(\land V\otimes \land (sV),D)$ with $D (sv) =
-s(dv)$. Furthermore a model for the equivariant free loop space
$LM\times_{S^1}ES^1$ is given by the commutative differential
graded algebra $$(\land V \otimes \land  sV \otimes \land u ,
D)\,,\quad D(u) = 0\,, \vert u\vert = 2\,, D(v) = d(v) + u s(v)\,,
D(sv) =-s(dv)\,.$$ A model for the inclusion $LM\to
LM\times_{S^1}ES^1$ is given by the projection $$\pi : (\land V
\otimes \land
 sV \otimes \land u , D) \to (\land V \otimes \land  sV , D)$$
defined by $\pi (u) = 0$. The map $s$ is the map induced in
cohomology by the   derivation   defined by $s(v) =sv$ and
$s(sv)=0$ for all $v\in V$. The Gysin sequence in cohomology is
then  identified   with a long exact sequence
 defined in terms of the Sullivan models
$$
\begin{array}{ccc}
H^n(LM) &\stackrel{\cong}{\longrightarrow} & H^n(\land (V\oplus
sV),D)\\ \uparrow {\scriptstyle H^n(p)}  & &\uparrow {\scriptstyle
H^n(\pi)} \\
 H^n( LM\times_{S^1}ES^1) &\stackrel{\cong}{\longrightarrow} &
 H^n(\land(V\oplus sV\oplus \mathbb Q u),D)\\
\uparrow {\scriptstyle c'}  & &\uparrow {\scriptstyle \cup u} \\
H^{n-2}(  LM\times_{S^1}ES^1)&\stackrel{\cong}{\longrightarrow} &
H^{n-2}(\land(V\oplus sV\oplus \mathbb Q u),D)\\ \uparrow
{\scriptstyle M'}  & &\uparrow {\scriptstyle s} \\ H^{n-1}(LM)
&\stackrel{\cong}{\longrightarrow} & H^{n-1}(\land (V\oplus sV),D)
\end{array}
$$

\vspace{3mm} Therefore we have proved

\vspace{3mm}\noindent {\bf Theorem 3.}  {\sl With rational
coefficients, the string bracket
 is dual
to the linear map $b^\vee :   H^*_{S^1}(LM) \to
H^*_{S^1}(LM)\otimes H^*_{S^1}(LM) $ defined by the composition $$
  H^*_{S^1}(LM) \stackrel{H^*(\pi)}{\longrightarrow} H^*(LM)
\stackrel{H^*(\mu_T\otimes1)\circ H^*(\Phi)^{-1}\circ
H^*(c)}{\longrightarrow}   H^*(LM)\otimes
H^*(LM)\stackrel{s\otimes s}{\longrightarrow} H^*_{S^1}(LM)\otimes
H^*_{S^1}(LM)\,. $$}

\vspace{5mm}\noindent {\bf Example.} Let $M$ be the product of two
   spheres of odd dimension $N$.
Models for $M$, $LM$ and $LM\times_{S^1}ES^1$ are then
  given by
$$
\begin{array}{ll}
M \, : & (\land (x,y),0)\,,\\ LM\,: & (\land (x,y, \bar x,\bar y),
0)\,,\\ LM\times_{S^1}ES^1\,: & (\land (x,y,\bar x,\bar y,u),
D)\,, Du = D(\bar x) =   D(\bar y) = 0\,, Dy=u\bar y\,, Dx = u\bar
x\,.
\end{array}
$$ A set of cocycles representing a basis of the vector space
$\tilde H_{S^1}^*(LM):= H^*_{S^1}(LM)/\mathbb Q u$ is composed of
the elements (\cite{BV}) $$
\renewcommand{\arraystretch}{1.5}
\left\{ \begin{array}{l} e_{a,b} = \bar x^a \bar y^b\,, \quad
(a,b)\in \mathbb N^2-(0,0) \\ f_{a,b} = (y\bar x - x\bar y)\bar
x^a\bar y^b\,, \quad (a,b) \in \mathbb N^2\,. \end{array} \right.
\renewcommand{\arraystretch}{1}$$ The vector space ${\tilde H}^*(LM)$ has the following basis $$
\left\{\begin{array}{l}
 e_{a,b}\\ f_{a,b}\end{array}\right.\,,\quad
 \left\{\begin{array}{l}   e'_{a,b} = xy\bar x^a\bar y^b,\,\,(a,b) \in \mathbb
N^2,\\ f'_{a,b} = x\bar x^a\bar y^b, \,\,
 (a,b) \in \mathbb N^2\\  f''_b
= y\bar y^b, \,\, b\in \mathbb N\end{array} \right.\,.
 $$
From the above description of the Gysin sequence we deduce $H^*(p)
(e_{a,b}) = e_{a,b}$, $  H^*(p)(f_{a,b}) = f_{a,b}$, $M'(f'_{a,b})
= e_{a+1,b}$,   $M'(e'_{a,b}) = f_{a,b}$, $M'(f''_b) = e_{0,b+1}$,
$M'(e_{a,b}) = 0$,  and $M'(f_{a,b}) = 0.$

 To fix signs, denote by $xy$ the fundamental class of $M$.
A straightforward computation shows that  $$
\renewcommand{\arraystretch}{1.9}
\begin{array}{l}
b^\vee (u^r) = 0\\
 b^\vee (\bar x^p\bar y^q) = \displaystyle\sum_{r=0}^p\sum_{s=0}^q
  \left(
\renewcommand{\arraystretch}{0.6}
\hspace{-1mm}\begin{array}{l} {\scriptstyle  p}\\{\scriptstyle
r}\end{array}\hspace{-1mm}
 \renewcommand{\arraystretch}{1}
\right)   \,
 \left(
\renewcommand{\arraystretch}{0.6}
\hspace{-1mm}\begin{array}{l} {\scriptstyle  q}\\{\scriptstyle
s}\end{array}\hspace{-1mm}
 \renewcommand{\arraystretch}{1}
\right)   \,\, \left( \bar x^r\bar y^{s+1} \otimes \bar
x^{p-r+1}\bar y^{q-s} - \bar x^{r+1}\bar y^s\otimes \bar
x^{p-r}\bar y^{q-s+1}\right)\\ b^\vee (f_{p,q}) = (f_{0,0}\otimes
1 + 1 \otimes f_{0,0})(\bar x\otimes \bar y -\bar y\otimes \bar
x)\cdot  b^\vee (\bar x^p\bar y^q)
\end{array}
\renewcommand{\arraystretch}{1}$$

To describe the string bracket in ${\cal H}_*(M)$ we choose dual
basis $t_r$, $a_{p,q}$ and $b_{p,q}$ of $u^r$, $\bar x^p\bar y^q$
and $f_{p,q}$.   In that basis the string bracket satisfies $$
\renewcommand{\arraystretch}{1.7}
\begin{array}{l}
\mbox{}  [b_{k,t}, a_{l,m}] =
 \left(
\renewcommand{\arraystretch}{0.6}
\hspace{-1mm}\begin{array}{c} {\scriptstyle  k+l}\\{\scriptstyle
k}\end{array}\hspace{-1mm}
 \renewcommand{\arraystretch}{1}
\right)   \,
 \left(
\renewcommand{\arraystretch}{0.6}
\hspace{-1mm}\begin{array}{c} {\scriptstyle  m+t}\\{\scriptstyle
t}\end{array}\hspace{-1mm}
 \renewcommand{\arraystretch}{1}
\right)   \, \frac{km-lt}{(k+l)(t+m)}\,\,\, b_{k+l-1,t+m-1}\,,\\
\mbox{}  [a_{k,t}, a_{l,m}] =
 \left(
\renewcommand{\arraystretch}{0.6}
\hspace{-1mm}\begin{array}{c} {\scriptstyle k+l}\\{\scriptstyle
k}\end{array}\hspace{-1mm}
 \renewcommand{\arraystretch}{1}
\right)   \,
 \left(
\renewcommand{\arraystretch}{0.6}
\hspace{-1mm}\begin{array}{c} {\scriptstyle  m+t}\\{\scriptstyle
t}\end{array}\hspace{-1mm}
 \renewcommand{\arraystretch}{1}
\right)   \, \frac{lt-km}{(k+l)(t+m)}\,\,\, a_{k+l-1,t+m-1}\,,\\
\mbox{} [b_{r,s}, b_{m,n}] = 0
 \end{array}
 \renewcommand{\arraystretch}{1}
 $$
In particular the string Lie algebra ${\cal H}_*(M)$ is not
nilpotent, since for instance
 $[a_{1,1}, a_{r,s}] = (r-s)\, a_{r,s}\,.$

\section{Hochschild cohomology of a Lie algebra}

 Let $(L,d_L)$ be a differential graded   Lie algebra.
The chain map $$\varphi : C_*L \to \bB (UL)$$ defined by
$$\varphi(sx_1\land \cdots \land sx_n) = \sum_{\sigma \in
\Sigma_n} \varepsilon_\sigma [x_{\sigma (1)}\vert \cdots \vert
x_{\sigma (n)}]\,,$$ is a quasi-isomorphism of coalgebras
(\cite{FHT}, Proposition 22.7), that can be extended for each
right $L$-module $P$ and each left $L$-module $N$ into a
quasi-isomorphism of complexes $$\Phi = 1 \otimes \varphi\otimes 1
: C_*(P;L;N) \to \bB (P, UL, N)\,.$$ In particular when $P = N =
UL$, $\Phi$ is a quasi-isomorphism of $UL$-bimodules.

By composition with $\Phi$ we get a quasi-isomorphism of complexes
$$\mbox{Hom}(\Phi, UL) : \mbox{Hom}_{(UL)^e}(\bB(UL;UL;UL), UL)
\to \mbox{Hom}_{(UL)^e}(C_*(UL;L;UL), UL)\,.$$ The linear
isomorphism $\mbox{Hom}_{(UL)^e}(C_*(UL;L;UL), UL) =
\mbox{Hom}_{\mathbb Q}(C_*L, UL)$ gives to the right hand side
complex   an algebra structure: $$f \cdot g : C_*L
\stackrel{\Delta}{\to} C_*L \otimes C_*L \stackrel{f\otimes
g}{\longrightarrow} UL\otimes UL \stackrel{\mu}{\to} UL\,.$$ The
left hand side complex is equipped with the Gerstenhaber usual
product.

\vspace{3mm}\noindent {\bf Theorem 4.} {\sl The correspondence
$\mbox{Hom}(\Phi, UL)$ is a quasi-isomorphism of differential
graded algebras inducing in homology an isomorphism of graded
algebras $$H\!H^*(UL,UL) \cong
H^*\left(\mbox{Hom}_{(UL)^e}(C_*(UL;L;UL), UL)\right)\,.$$ }

\vspace{3mm}\noindent {\bf Proof.} We form the commutative diagram
$$
\renewcommand{\arraystretch}{1.5} \begin{array}{ccc}
\mbox{Hom}_{(UL)^e}(\bB(UL;UL;UL), UL) &\stackrel{ \scriptstyle
 {\scriptsize  Hom}(\Phi, UL)}{\longrightarrow} &
\mbox{Hom}_{(UL)^e}(C_*(UL;L;UL), UL)\\ \downarrow \cong &&
\downarrow \cong\\ \mbox{Hom} (\bB(UL), UL) &\stackrel{
{\scriptsize Hom}
 \scriptstyle
(\varphi, UL)}{\longrightarrow} & \mbox{Hom} (C_*(
L ), UL)\,. \end{array}
\renewcommand{\arraystretch}{1} $$ Since $\varphi$ is a morphism
of coalgebras, $\mbox{Hom}(\varphi, UL)$ is a morphism of
algebras.
 \hfill $\square$

 \vspace{3mm} The canonical isomorphism
 $$\mbox{Hom}_{(UL)^e}(C_*(UL;L;UL), UL) \cong \mbox{Hom}_{\mathbb Q}
 ( C_*(UL;L;UL)\otimes_{(UL)^e} (UL)^\vee, \mathbb Q)$$
 gives then clearly the isomorphism
 $$H\!H^*(UL;UL) \cong \mbox{Hom}\left(H_*\left(C_*(UL;L;UL)\otimes_{(UL)^e} (UL)^\vee\right),\mathbb Q\right)\,.$$
 We equip   $UL$ and  $(UL)^\vee$ with the left adjoint representation, $$l\cdot x = [l,x]\,,
  (l\cdot f)(x) = - (-1)^{\vert l\vert \cdot \vert f\vert} f([l,x])\,,\quad
  l\in L, x\in UL, f\in (UL)^\vee\,,$$
 and, to avoid confusion,  we denote these modules respectively by $(UL)_a$ and $(UL)_a^\vee$.
 Then a straightforward computation shows

 \vspace{3mm}\noindent {\bf Lemma 1.}  {\sl The linear isomorphism
 $$C_*(L, (UL)_a^\vee) = C_*(UL;L;UL)\otimes_{(UL)^e} (UL)^\vee$$
 is an isomorphism of chain complexes.}

 \vspace{3mm}  The chain complex  $C_*(L;(UL)_a^\vee)$ is
a coalgebra that is the tensor product of the coalgebras
 $C_*(L)$ and $(UL)^\vee$. The dual  algebra is the cochain algebra $C^*(L;(UL)^\vee_a)$. Since
 $\mbox{Hom}(C_*(L;(UL)^\vee_a),\mathbb Q) =
 \mbox{Hom}_{UL}(C_*(L;UL), (UL)_a) = C^*(L;(UL)_a)$,
  we deduce

 \vspace{3mm}\noindent {\bf Corollary 1.}  {\sl Let $L$ be a
 differential graded Lie algebra, then we have   isomorphisms of
 algebras
 $$H\!H^*(UL,UL) \cong \mbox{Hom} \left(H_*( C_*(L;(UL)_a^\vee)),
 \mathbb Q\right)\cong H^*(C^*(L;(UL)_a))\,.$$}

 \section{Lie models and the cap product $C^*(L;M) \to C_{m-*}(L;M)$}

\subsection{Minimal Lie models}

Let $S$ be a simply connected space. By Quillen theory
(\cite{Q},\cite{Maj},\cite{FHT}), we can associate to $S$ a
differential graded Lie algebra $L = {\cal L}_S$, that is  free as
a graded Lie algebra, $L = \mathbb L(V)$, with a   differential
$d$ satisfying $d(V) \subset \mathbb L^{\geq 2}(V)$ (the sub
vector space generated by the brackets of length at least two).

The differential graded Lie algebra  $L$ contains all the rational
homotopy type of $S$. In particular there exists a
quasi-isomorphism from the enveloping algebra on $L$ into the
chain algebra $C_*(\Omega S;\mathbb Q)$, $$\varphi_S : UL
\stackrel{\simeq}{\to} C_*(\Omega S;\mathbb Q)\,.$$ Moreover a
differential graded Lie algebra $L$ is a Lie model for $S$ if and
only if   the cochain algebra on $L$, $ C^*L$, is a Sullivan
  model of $S$.

\subsection{A Lie model for $LM$}

Starting from a differential graded Lie algebra $L$ we construct
another differential graded Lie algebra $L^S$, that is a Lie model
for $LM$ when $L$ is a Lie model for $M$.

The differential graded Lie algebra $  L^S$ is defined as follows
$$L^S_n = L_n \oplus \overline{L}_n\,, \quad \overline{L}_n =
L_{n+1}\,,\quad d\overline{x} = - \overline{dx}\,, \quad
(-1)^{\vert a\vert} [a,\overline{b}] = \overline{[a,b]}\,,\quad
[\overline{a},\overline{b}] = 0\,.$$ In particular $\overline{L}$
is an abelian sub Lie algebra. By the Poincar\'e-Birkhoff-Witt
Theorem the chain coalgebra $C_*\overline{L} = \land
(s\overline{L})$ is isomorphic to $UL$ as a differential graded
coalgebra: $$\varphi : C_*(\overline{L})
\stackrel{\cong}{\longrightarrow} UL\,, \quad \varphi (s\bar
x_1\land \cdots \land s\bar x_n) = x_1 \cdots x_n\,.$$ Since
$C_*(L^S) \cong C_*(L)\otimes C_*(\overline{L})$, we   have by
direct computation

\vspace{3mm}\noindent {\bf Lemma 2.} {\sl The morphism $1\otimes
\varphi : C_*(L^S)   \to C_*(L;(UL)_a)$ is an isomorphism of
coalgebras.}

\vspace{3mm}The relation with the free loop space is the content
of  Theorem 5.

\vspace{3mm}\noindent {\bf Theorem 5.} {\sl If $L$ is a Lie model
for $M$, then $L^S$ is a Lie model for $LM$, and a Sullivan
model for $LM$ is given by the commutative differential graded
algebra $C^*(L;(UL)^\vee_a)$.}

\vspace{3mm}\noindent {\bf Proof.}  Since $L$ is a Lie model for $M$,
  $C^*(L) =(\land V,d)$, with $V =
s(L^\vee)$ (\cite{FHT}, p. 320), is a Sullivan model for $M$. We
remark that $$C^*(L^S) = (\land (V\oplus sV),D)$$ with $D(v) =
d(v)$ and $D(sv) = -sd(v)$. Therefore $L^S$ is a Lie model for
$LM$.

The last assertion follows from the following sequence of
isomorphisms of differential graded algebras $$
\renewcommand{\arraystretch}{1.5}
\begin{array}{ll}
C^*(L^S) \cong \mbox{Hom}(C_*(L;(UL)_a),\mathbb Q) &=
\mbox{Hom}(C_*(L;UL) \otimes_{UL}(UL)_a,\mathbb Q) \\ & =
\mbox{Hom}_{UL}(C_*(L;UL), (UL)_a^\vee)= C^*(L;(UL)_a^\vee)\,.
\end{array}
\renewcommand{\arraystretch}{1}$$

\hfill $\square$

\subsection{A cochain model for the path composition $LM\times_MLM\to LM$}

The multiplication $\mu : UL\otimes UL\to UL$  is a morphism of
$UL$-modules for the adjoint representation. Here $UL$ acts
diagonally by adjunction on $UL\otimes UL$. Therefore    we have a
morphism of cochain complexes
$$C^*(L; \mu^\vee):  C^*(L;(UL)_a^\vee)
\to C^*(L;(UL)_a^\vee \otimes (UL)^\vee_a)\,.$$

\vspace{3mm}\noindent {\bf Theorem 6.}  {\sl When $L$ is a Lie
model for $M$, then the morphism $C^*(L; \mu^\vee) $ is a morphism
of commutative differential graded algebras and is a Sullivan
model for the path composition $LM\times_MLM \to LM$.}

\vspace{3mm}\noindent {\bf Proof.} The complex $C_*({\mathbb
Q};L;UL)$ is a differential graded coalgebra whose coproduct is
the tensor product of the coproduct  on $C_*L$ and the coproduct
on $UL$.  The diagonal $\Delta$ on $C_*L$ induces a
quasi-isomorphism of differential graded coalgebras $$\Delta' :
C_*L \to C_*({\mathbb Q};L;UL) \otimes_{UL} C_*(UL;L;\mathbb
Q)\,$$
 defined   by $\Delta'(x)=\sum_i x_i\otimes 1 \otimes x_i'$ when $\Delta (x) =
 \sum_i x_i\otimes x_i'$.

The dual map $(\Delta')^\vee : C^*L\otimes (UL)^\vee\otimes C^*L
\to C^*L$ is a quasi-isomorphism of differential graded algebras,
and the injection $1 \otimes\varepsilon'\otimes 1 : C^*L\otimes
C^*L \to C^*L\otimes (UL)^\vee\otimes C^*L$ is a relative Sullivan
model for the product $\mu : C^*L\otimes C^*L \to C^*L$. Here
$\varepsilon  : UL\to \mathbb Q  $ denotes  the augmentation of
$UL$ and $\varepsilon'$ the dual map.

 A cochain model for the diagram
$$
\begin{array}{ccccc}
M & \stackrel{\Delta}{\longrightarrow} & M \times M & \stackrel{\Delta}{\longleftarrow} & M\\
\parallel && \uparrow{\scriptstyle \rho} &&\parallel\\
M &\stackrel{\Delta}{\longrightarrow} & M \times M\times M & \stackrel{\Delta}{\longleftarrow}
& M
\end{array}
$$
 is given by
$$
\renewcommand{\arraystretch}{1.5}
\begin{array}{ccccc}
 C^*L & \stackrel{\mu}{\leftarrow} & C^*L\otimes C^*L & \stackrel{\mu}{\to} & C^*L  \\
  \parallel && {  \scriptstyle 1 \otimes \tau\otimes 1  }\downarrow  && \parallel  \\
 C^*L & \stackrel{(\mu\otimes 1)\circ \mu}{\leftarrow}&   C^*L \otimes C^*L\otimes
 C^*L
& \stackrel{(\mu\otimes 1)\circ \mu}{\to} & C^*L\,.
\end{array}
\renewcommand{\arraystretch}{1}
$$ Here     $\rho (a,b,c)= (a,c)$  and  $\tau$  denotes the unit
$\mathbb Q \to C^*L$. In order to replace the maps $\mu$ and
$(\mu\otimes 1)\circ \mu$  by relative models, we decompose the
right hand square into two squares $$
\renewcommand{\arraystretch}{1.5}
\begin{array}{ccccc}
  C^*L\otimes C^*L  & \stackrel{  1 \otimes \varepsilon'\otimes 1}{\to}
  & C^*L\otimes (UL)^\vee\otimes C^*L   &\stackrel{(\Delta')^\vee}{\to} &C^*L \\
  {\scriptstyle 1 \otimes \tau\otimes 1}\downarrow   &&
\mbox{} \hspace{15mm} \downarrow {\scriptstyle  1 \otimes  \nabla'\otimes 1  }&&\parallel\\
  C^*L \otimes C^*L\otimes C^*L &
  \stackrel{\bar\varepsilon}{\to} & C^*L\otimes (UL)^\vee\otimes C^*L \otimes (UL)^\vee \otimes
  C^*L
&\stackrel{((\Delta')^\vee\otimes 1)(\Delta')^\vee}{\to} &C^*L
\end{array}
\renewcommand{\arraystretch}{1}
$$
where $C^*L\otimes (UL)^\vee\otimes C^*L =
(C_*(L;UL)\otimes_{UL}C_*(UL;L))^\vee$, $\bar\varepsilon = 1
\otimes \varepsilon'\otimes 1\otimes \varepsilon'\otimes 1$,
$C^*L\otimes (UL)^\vee\otimes C^*L \otimes (UL)^\vee \otimes C^*L
= \left(C_*(L;UL)\otimes_{UL}
C_*(UL;L;UL)\otimes_{UL}C_*(UL;L)\right)^\vee$, $\nabla$ is the
coproduct on $(UL)^\vee$ and $\nabla'$ is the composition
$$\nabla' : (UL)^\vee \stackrel{\nabla}{\longrightarrow} (UL)^\vee\otimes \mathbb Q\otimes
(UL)^\vee  \stackrel{1\otimes \tau\otimes 1}{\longrightarrow}
(UL)^\vee\otimes C^*L \otimes (UL)^\vee = (C_*(UL;L;UL))^\vee
  \,.$$

Now we consider the following commutative diagram in which the
front face and the back face  are fiber squares. By section 3.5,
this implies that $C^*(L;\mu^\vee)$ is a model for the path
composition $LM\times_MLM \to LM$.

\mbox{}\hspace{15mm}\begin{pspicture}(0,0)(9,6)
\psline{->}(7.6,0)(3.5,0)
\psline{->}(8.4,3)(4.5,3)\psline{->}(3.9,5)(2,5)\psline{->}(4.9,2)(1.5,2)
\psline{->}(3,0.4)(3,2.7)
\psline{->}(1,2.4)(1,4.7)\psline{->}(6,2.4)(6,4.7)\psline{->}(1.3,1.8)(2.7,0.2)
\psline{->}(1.3,4.8)(2.7,3.2)\psline{->}(6.3,4.8)(8.7,3.2)\psline{->}(6.3,1.8)(8.7,0.2)
\psline{->}(9,0.4)(9,2.6) \rput(6,5){$\scriptstyle
(C_*(L;UL)\otimes_{UL}C_*(UL;L))^\vee$} \rput(1,5){$\scriptstyle
C^*(L;(UL)_a)$} \rput(1,2){$C^*L$} \rput(6,2){$\scriptstyle
C^*L\otimes C^*L$}
 \rput(6.7,3.3){$\scriptstyle 1\otimes \varepsilon'\otimes 1$}\rput(4,2.2){$\scriptstyle \mu$}
\rput(3,3){$\scriptstyle C^*(L;(UL)_a\otimes (UL)_a)$}
\rput(3,0){$C^*L$} \rput(9,0){$\scriptstyle C^*L\otimes
C^*L\otimes C^*L$} \rput(9,3){$\scriptstyle E$}
\rput(5,0.3){$\scriptstyle(\mu\otimes 1)\mu$}
\rput(2.4,1.4){$\scriptstyle id_{C^*L}$}
\rput(9.2,2.5){$\scriptstyle \bar\varepsilon$}
\rput(7.9,1.2){$\scriptstyle 1\otimes \tau\otimes 1
$}\rput(7.9,4.3){$\scriptstyle 1\otimes \nabla'\otimes 1$}
\rput(2.8,4){$\scriptstyle C^*(L;\mu^\vee)$}
\end{pspicture}

\vspace{4mm}\noindent with $E = \left(C_*(L;UL)\otimes_{UL}
C_*(UL;L;UL)\otimes_{UL}C_*(UL;L)\right)^\vee$. \hfill$\square$

\vspace{4mm}  Denote by $ L_1^S$ and $ L_2^S$
two copies of $  L^S$, and  $ L^T =   L_1^S\oplus_{L}
  L_2^S = L \oplus \bar L_1 \oplus \bar L_2$. We denote by
$\pi :  L^T \to   L^S$ the projection obtained by mapping
identically each $\ L_i^S $ to $ L^S$. We observe that $C_*( L^T)$
is isomorphic to $C_*(L;(UL)_a^\vee \otimes (UL)_a^\vee)$ and that
the following diagram commutes $$\renewcommand{\arraystretch}{1.4}
\begin{array}{ccc}
C^*(L;(UL)_a^\vee) & \stackrel{C^*(L;\mu^\vee)}{\longrightarrow} &
C^*(L;(UL)_a^\vee \otimes (UL)^\vee_a)\\
\parallel &&\parallel\\
\mbox{Hom}(C_*(L;(UL)_a), \mathbb Q) &
\stackrel{\mbox{Hom}(C_*(L;\mu),\mathbb Q)}{\longrightarrow} &
\mbox{Hom}(C_*(L;(UL)_a\otimes(UL)_a),\mathbb Q)\\
\parallel && \parallel \\
C^*(L^S) & \stackrel{C^*(\pi)}{\longrightarrow} & C^*(L^T)
\end{array}
\renewcommand{\arraystretch}{1}
$$
This shows that $\pi : L^T \to L^S$ is a Lie model for the path composition
$LM\times_MLM\to LM$.

\vspace{3mm} Recall that a coformal space is a space that admits a
Sullivan minimal model with a purely quadratic differential. The
minimal  Sullivan model of a coformal space is the cochain algebra
on the Lie algebra $\pi_*(\Omega M)\otimes \mathbb Q$. Therefore,

\vspace{2mm}\noindent {\bf Corollary  2.}  {\sl Let $M$ be a
coformal manifold with
 minimal model   $(\land V,d)$, then a model for the path composition
$LM\times_MLM \to LM$ is given by $$id\otimes \Delta : (\land  V
\otimes \land  \bar V,D ) \to (\land  V\otimes \land \bar V
\otimes \land \bar V'),D)\,,$$ where $\Delta$ is the coproduct in
$\land \bar V=H^*(\Omega M)$ and $(\land V \otimes \land \bar
V,D)$ is the model of the free
 loop space.}

 \vspace{3mm}\noindent {\bf Example.}
Let $M$ be the $11$-dimensional manifold   obtained by taking the
pullback of the tangent sphere bundle to $S^6$ along the map $f :
S^3 \times S^3 \to S^6$ that collapses the 3-skeleton into a
point. $$
\begin{array}{ccc}
M & \to &\tau_SS^6\\ \downarrow & &\downarrow \\ S^3\times S^3
&\stackrel{f}{\to} & S^6
\end{array}
$$ The minimal model of $M$ is
 $$(\land (x,y,z), d)\,, \quad dx= dy = 0 \,,  dt = xz\,,\hspace{5mm} \vert x\vert =\vert y\vert = 3\,,
 \vert z\vert = 5\,.$$ Thus $M$ is a coformal space.
  A model for the path composition
 $M^{[0,1]}\times_MM^{[0,1]} \to M^{[0,1]}$
is given by
 $$\varphi : (\land (x,y,z,x',y',z',\bar x,\bar y,\bar z),D) \to
 (\land (x,y,z,x',y',z',x'',y'',z'',
 \bar x,\bar y,\bar z,\bar x',\bar y',\bar z'),D)\,, $$
 with
 $D(\bar x) = x-x'$, $D(\bar y) = y-y'$,   $D(\bar z) = z-z' -\frac{1}{2} \bar x (y+y')
  +\frac{1}{2} (x+x')\bar y$,
   $D(\bar x') = x'-x''$, $D(\bar y') = y'-y''$,   $D(\bar z') = z'-z'' -\frac{1}{2} \bar x'
    (y'+y'')
  +\frac{1}{2} (x'+x'')\bar y'$,
 $\varphi (\bar x) = \bar x+ \bar x'$, $\varphi (\bar y)= \bar y+ \bar y'$, $\varphi (\bar z) =
 \bar z+ \bar z' + \frac{1}{2}\bar x\bar y'-\frac{1}{2}\bar x'\bar y$.
   The induced model for the path composition
 $LM\times_MLM \to LM$
is then given by
 $$\varphi : (\land (x,y,z, \bar x,\bar y,\bar z),D) \to
 (\land (x,y,z,
 \bar x,\bar y,\bar z,\bar x',\bar y',\bar z'),D)\,, $$
 with
 $D(\bar x) = 0$, $D(\bar y) = 0$,   $D(\bar z) =   -  \bar x  y
  +   x \bar y$,
   $D(\bar x') = 0$, $D(\bar y') = 0$,   $D(\bar z') =   -  \bar x'
     y
  + x \bar y'$,
 $\varphi (\bar x) = \bar x+ \bar x'$, $\varphi (\bar y)= \bar y+ \bar y'$, $\varphi (\bar z) =
 \bar z+ \bar z' + \frac{1}{2}\bar x\bar y'-\frac{1}{2}\bar x'\bar y$.

\subsection{The cap product}

Let $L$ be a differential graded Lie algebra and $N$ be a left
$UL$-module. Then $N$ is a right $UL$-module for the action
defined by $$n\cdot x := -(-1)^{\vert n\vert \cdot \vert x\vert}
x\cdot n\,.$$ This gives the opportunity to consider in the same time a module   as a
left module and   as a right module.

Let $c = \sum_i sx_{i_1}\land \cdots \land sx_{i_q} $ be a cycle
of degree $m$ in $C_qL$. We define the cap product by $c$, $$
cap_c : C^{q-r}(L;N) = \mbox{Hom}_{UL}(C_{q-r}(L;UL),N) \to C_r(L;
N)$$ by mapping an element $f$ to the element $f\cap c$ defined by
$$ f\cap c = (-1)^m\,\,\sum_i\sum_{\sigma\in \Sigma_q}   (-1)^{\vert f\vert
\cdot (\vert sx_{\sigma(i_1)} \land \cdots \land sx_{\sigma
(i_r)\vert})}\, \varepsilon_\sigma\,\, sx_{\sigma(i_1)} \land
\cdots \land sx_{\sigma (i_r)} \otimes f(sx_{\sigma(i_{r+1})}
  \cdots   sx_{\sigma (i_q)})\,.$$

Remark that in $ C^*(L; N)$,  $N$ is considered as a right
$UL$-module and in $C_*(L;N)$,  $N$ is viewed as a left
$UL$-module. By a standard computation  we then have

\vspace{3mm}\noindent {\bf Lemma 3.} {\sl The morphism $ cap_c$ is
a natural morphism of complexes: If $g : P \to N$ be a morphism of
left $UL$-modules, then the following diagram is commutative $$
\renewcommand{\arraystretch}{1.4}
\begin{array}{ccc}
C_*(L;P)  & \stackrel{ C_*(L;g)}{\longrightarrow} &  C_*(L;N)\\
{\scriptstyle  cap_c}{\uparrow} & &   {\scriptstyle  cap_c}{\uparrow} \\
 C^*(L; P) & \stackrel{C^*(L;g)}{\longrightarrow} &  C^*(L;N)
\end{array}
\renewcommand{\arraystretch}{1}
$$
}

\vspace{3mm}

In particular since the multiplication $\mu : UL\otimes UL\to UL$,
and its dual co-multiplication $\mu^\vee$
 are morphisms of $UL$-modules
for the adjoint representation, we have

\vspace{3mm}\noindent {\bf Theorem 7.}  {\sl The following diagram
commutes $$
\renewcommand{\arraystretch}{1.4}
\begin{array}{ccc}
C_*(L;(UL)_a^\vee)  & \stackrel{
C_*(L;\mu^\vee)}{\longrightarrow} &
 C_*(L;(UL)_a^\vee\otimes (UL)_a^\vee)\\
{\scriptstyle  cap_c}{\uparrow} & &   {\scriptstyle  cap_c}{\uparrow} \\
\mbox{Hom}_{UL}(C_*(L;UL),(UL)_a^\vee) &
\stackrel{\mbox{\scriptsize Hom}(1,\mu^\vee)}{\longrightarrow}
& \mbox{Hom}_{UL}(C_*(L;UL),(UL)_a^\vee\otimes (UL)_a^\vee) \\
\parallel & & \parallel\\
C^*(L;(UL)_a^\vee) & \stackrel{C^*(L; \mu^\vee)}{\longrightarrow} &
C^*(L;(UL)_a^\vee \otimes (UL)^\vee_a)
\end{array}
\renewcommand{\arraystretch}{1}
$$}

\section{Proof of Theorem B}

The proof is based on the following diagram
$$
\renewcommand{\arraystretch}{1.6}
\begin{array}{ccccc}
{\scriptstyle C_*(L;(UL)^\vee_a )} &
\stackrel{C_*(L;\mu^\vee)}{\longrightarrow}
 & {\scriptstyle  C_*(L;(UL\otimes UL)^\vee_a)} &
 \stackrel{C_*(\Delta;id)}{\longrightarrow}
& {\scriptstyle  C_*(L;(UL)^\vee_a)\otimes C_*(L;(UL)^\vee_a)}
\\
\uparrow {\scriptstyle  cap_c} && \uparrow {\scriptstyle  cap_c} &&
\\
{\scriptstyle  C^*(L;(UL)_a^\vee)} &
\stackrel{C^*(L;\mu^\vee)}{\longrightarrow} & {\scriptstyle
C^*(L;(UL\otimes UL)^\vee_a)} &
&\uparrow {\scriptstyle ( cap_c)\otimes (  cap_c)}\\
&& \uparrow{\scriptstyle \cong} &&\\
&&{\scriptstyle  C^*(L;(UL)^\vee_a )\otimes C^*(L;(UL)^\vee_a
)\otimes \land \bar V } &\stackrel{\mu_T\otimes 1}{\to} &
{\scriptstyle C^*(L;(UL)_a^\vee)\otimes C^*(L;(UL)_a^\vee)}\,,
\end{array}
\renewcommand{\arraystretch}{1}
$$
 where $c$ is a cycle representing the orientation class in homology.
The dual of the upper row induces    the Gerstenhaber product in
Hochschild cohomology (Section 5, Corollary 1), the left hand
square is commutative by Theorem 7, and the bottom line induces in
cohomology the dual of the Chas-Sullivan product
 (Theorem A and Theorem 6). Since $H^*(L;\mathbb Q)$ satisfies Poincar\'e duality, a classical
  spectral sequence argument shows that $cap_c$ is a
  quasi-isomorphism.
 We will prove in this chapter (Theorem 8) that the right hand square commutes in
homology. This will end  the proof of Theorem A. \hfill $\square$

\subsection{The cap product and the diagonal class}

Let $L$ be a Lie model for $M$. The algebra
$H^*(C^*(L))=H^*(M;\mathbb Q)$  is a Poincar\'e duality algebra
with a fundamental class in degree $m$. According to the notations
of section 3.4, we denote by $u =u_M$ a cycle in $C_*(L)$ whose
class is a generator of $H_m(C_*L)$ and by $\alpha_i$, $i=1,
\ldots  , p$,  a graded basis of $H^*(L;\mathbb Q)$. The
Poincar\'e dual basis $\alpha_i^{\#}$ is defined by $$\langle
\alpha_i\cup\alpha_j^{\#}; [u_M]\rangle = \left\{
\begin{array}{ll} 1 & \mbox{if }\, i=j\\ 0 & \mbox{if}\, i\neq
j\end{array} \hspace{1cm}  \mbox{(\cite{MS}, Theorem 11.10)}
\right.$$ We denote by $a_i$ and $a_i^{\#}$ cocycles representing
the classes $\alpha_i$ and $\alpha_i^{\#}$. The diagonal
cohomology class $\delta_M$ is  then the class of the cocycle
$T=\sum_i (-1)^{\vert a_i\vert} a_i\otimes a_i^{\#}$ (\cite{MS} ,
Theorem 11.11).

\vspace{3mm} The purpose of this section consists of the proof of
the following theorem.

\vspace{3mm}\noindent {\bf Theorem 8.} {\sl With the above
notations, the following diagram commutes in homology $$
\begin{array}{ccc}
C_*(L;(UL)^\vee_a\otimes (UL)^\vee_a) &
\stackrel{C_*(\Delta;id)}{\longrightarrow} & C_*(L;(UL)^\vee_a)
\otimes C_*(L;(UL)^\vee_a)\\ \uparrow {\scriptstyle  ( cap_{u })}
&&\\ C^*(L;(UL)^\vee_a\otimes (UL)^\vee_a) && \uparrow
{\scriptstyle ( cap_{u })\otimes ( cap_{u })}\\
{\scriptstyle\Phi}\uparrow{\scriptstyle \cong} &&\\
C^*(L;(UL)^\vee_a)\otimes C^*(L;(UL)_a^\vee)\otimes \land \bar V &
\stackrel{\mu_T\otimes 1}{\longrightarrow}&
C^*(L;(UL)^\vee_a)\otimes C^*(L;(UL)^\vee_a)
\end{array}
$$
The morphisms $\mu_T$ and $\Phi$ are defined in section 3.4. }

\vspace{3mm} For the proof we will use the following  general
Lemma  for   modules over  differential graded algebras.

\vspace{3mm}\noindent {\bf Lemma 4.}  {\sl Let $R$ be a
differential graded algebra, $R\otimes V$ be a semifree
$R$-module, and $f,g : N \to P$   morphisms of right $R$-modules.
If $H_*(f) = H_*(g)$ then the maps $f\otimes 1, g\otimes 1 :
N\otimes_R (R\otimes V) \to P\otimes_R (R \otimes V)$ induce also
the same map in homology.}

\vspace{3mm}
We also recall the relation between the diagonal class and the cap product
with the fundamental class (the proof is a direct computation following the previous definitions).

\vspace{3mm}\noindent {\bf Lemma 5.} {\sl With the previous
notations the following diagram is commutative
$$\begin{array}{ccc} H_*( L;\mathbb Q)& \stackrel{H_*(\Delta
)}{\longrightarrow} & H_*(L\oplus L;\mathbb Q)\\
\uparrow{\scriptstyle  ( cap_{[{u }]})} && \uparrow {\scriptstyle
(  cap_{[{u }]})\otimes
 ( cap_{[u]}) }\\
 H^*(L;\mathbb Q) &\stackrel{{\mu_T}^*}{\to} & H^*(L;\mathbb Q)\otimes H^*(L;\mathbb Q)
 \end{array}
 $$}

\vspace{3mm}\noindent{\bf Proof of Theorem 8.}

We denote $R= C^*L\otimes C^*L$. We remark that   $R\otimes V =
C^*(L;(UL)^\vee_a)\otimes C^*(L;(UL)^\vee_a)$ is a  semifree
$R$-module. Remark also that   diagram in Theorem 8 is     the
tensor product of a diagram of $R$-module by $R\otimes V$ over
$R$.
  By Lemma 5,
 this diagram commutes in homology when $V =\mathbb Q$.
It commutes then in homology by Lemma 4, and this achieves the
proof of the Theorem. \hfill$\square$

\subsection{Proof of Lemma 4}

Let $R$ be a differential graded algebra. Each    differential
$R$-module $N$ admits a semifree resolution, i.e., there exists a
semifree $R$-module $S$ and a quasi-isomorphism $\psi : S \to N$.
For recall $S$ is a semifree resolution means that $S =
\sum_{n\geq 0}  R\otimes W(n)$ with $d(W(n))\subset \sum_{i<n}
R\otimes W(i)$. The sub-complexes $F^p = \sum_{q\leq p} R \otimes
W(q)$ form then a   filtration that induces a spectral sequence
converging to
 the cohomology of $N$,
with $E^0_{q,*} = R \otimes W(q)$.

\vspace{2mm} It is easy to see that each module $N$ admits a
semifree resolution $P$ such that the associated spectral sequence
satisfies $E^2_{0,*} = H_*(M)$ as an $H_*(R)$-module and
$E^2_{q,*} = 0$ for $q>0$.

\vspace{3mm}  We can assume that $R\otimes V$ is a  semifree
resolution satisfying this property. The
  spectral sequences obtained by filtering $N \otimes_R R\otimes V$ and $P\otimes_R R \otimes V$ by
   the filtration degree in $V$, $F^p$, degenerates at the $E_2$-term.
  Since $E_1(f) = H(f)\otimes id_V $ and $E_1(g) = H(g)\otimes id_V $, the   maps $f$ and $g$ induce the same
  map
at the  $E_1$-level.
   Since the spectral sequences degenerate at the $E_2$-terms
   and the morphisms $f$ and $g$ are the identity on $R\otimes V$, $H(f) = H(g)$.
   \hfill $\square$

\section{The intersection morphism}

For each subspace $Z \hookrightarrow M$, we denote by $L_ZM$ the space of loops in $LM$ that
originates in $Z$. We have thus the commutative pullback diagram
$$
\begin{array}{ccc}
L_ZM & \stackrel{i}{\longrightarrow} & LM\\
\downarrow && \downarrow\\
Z & \hookrightarrow & M\,.
\end{array}
$$
Denote now by $D$ a closed disk around the base point $x_0$. Since $D$ is contractible,
$L_DM$
is homotopy equivalent to $D \times \Omega M$. The restriction morphism
 $I : H_*(LM) \to H_{*-m}(\Omega M)$ is the composition
 $$\tilde H_*(LM) \to H_*(LM, L_{M-x_0}M) \stackrel{Excision}{\longrightarrow }
 H_*(L_DM,L_{S^{m-1}}M) \cong H_m(D, S^{m-1}) \otimes H_{*-m}(\Omega M)\,.$$

 Let $(\land V \otimes \land\bar V,D)$ be a  Sullivan
  model for the free loop space, and let
 $\omega \in (\land V)^m$ be a cocyle representing the fundamental class. Then the direct sum
 $J = (\land V)^{>m}  \oplus S$, where $(\land V)^m = \mathbb Q\cdot \omega\oplus S$,
 is an acyclic differential ideal in $\land V$; the quotient map $q : (\land V,d) \to (A,d)
  = (\land V / J,d)$ is a quasi-isomorphism.
  We form the tensor product
  $$(A\otimes \land \bar V,D) :=
  (A,d)\otimes_{\land V}(\land V\otimes \land \bar V,D) \,.$$

 \vspace{2mm}\noindent {\bf Lemma 6.}  {\sl The injection
 $i: \land \bar V\to A\otimes \land \bar V$
 defined by $i(a) = (-1)^m \omega \otimes a$ is a morphism of complexes of degree $m$
 inducing in cohomology the dual of the intersection morphism, $I^\vee :
 H^*(\Omega M) \to H^{*+m}(LM)$.}

 \vspace{3mm}\noindent {\bf Proof.} A model for the injection $L_{M-x_0}M \hookrightarrow LM$
 is given by the quotient map $q: (A \otimes \land \bar V,D)\to
 (A/(\omega) \otimes \land \bar V,D)$. Therefore a model for the cochain complexes
  $C^*(LM,L_{M-x_0}M)$ is given by Ker $q =(\mathbb Q\cdot \omega \otimes \land \bar V,0)$.
  This implies that a model for the composition
  $$H^*(\Omega M) \stackrel{\cong}{\longrightarrow} H^{*+m}(LM,L_{M-x_0}M) \to H^{*+m}(LM)$$
  is given by the composition
  $$\land \bar V \hspace{3mm}\stackrel{a\mapsto (-1)^m \omega\otimes a}{\longrightarrow}
   \hspace{4mm}
  \mathbb Q\cdot \omega\otimes\land \bar V \to H^*(A\otimes \land \bar V,D)\,.$$
  \hfill$\square$

   \vspace{3mm}\noindent {\bf Theorem D.}  {\sl There exists an isomorphism of algebras
 making commutative the diagram
 $$
 \begin{array}{ccc}
   H_*(LM) &\stackrel{ \Phi\,\cong}{\longrightarrow}& H\!H^*(C^*M,C^*M)\\
   {\scriptstyle I}\downarrow && \downarrow {\scriptstyle H\!H^*(C^*M,\varepsilon)}\\
   H_*(\Omega M) & \stackrel{\Psi}{\longrightarrow} & H\!H^*(C^*M;\mathbb Q)
   \end{array}$$
   where $\varepsilon : C^*M\to \mathbb Q$ denotes the usual augmentation.}

 \vspace{3mm}\noindent {\bf Proof.}

  Denote now by $L$ the minimal Lie model of $M$. The algebra
   $C^*(L;(UL)^\vee_a) = C^*(L)\otimes (UL)^\vee$ is then a
   free model for the free loop space $LM$, and we can apply Lemma 6 to this model to have a
   model for the restriction morphism. Now, since the projection $q : C^*(L;(UL)_a^\vee)
   \to A\otimes (UL)_a^\vee$ is a quasi-isomorphism, $q$ admits a lifting $q' : (UL)^\vee \to
   C^*(L;(UL)^\vee_a)$.

  Denoting by $c$ a fundamental class in homology with $\langle \omega,c\rangle = 1$,
   we obtain a commutative diagram
   $$\begin{array}{ccc}
   (UL)^\vee & \stackrel{q'}{\longrightarrow} & C^*(L;(UL)^\vee_a)\\
   \parallel && \downarrow {\scriptstyle (-1)^m -\cap c}\\
   (UL)^\vee & \stackrel{e}{\longrightarrow} & C_*(L;(UL)_a^\vee)\,,
   \end{array}
   $$
   where $e(a) = 1 \otimes a$.

   The dual diagram yields in homology to a diagram of algebras whose vertical maps
   are isomorphisms
   $$
   \begin{array}{cccc}
   H_{*}(\Omega M) &\stackrel{I}{\longleftarrow} & H_{*+m}(LM)\\
   \uparrow && \uparrow\\
   H_*(UL) & \stackrel{e'}{\longleftarrow} & H^*(\mbox{Hom}(C_*(L;(UL)^\vee_a), \mathbb Q))&=H^*(C^*(L;(UL)_a))\,.
   \end{array}\hspace{2mm}(*)
   $$

   Remark that the two $C_*(L)$-bimodules $ C_*(L;UL)\otimes_{UL}C_*(UL;L)$
   and
 $C_*L$ are quasi-isomorphic. Therefore we have by duality a quasi-isomorphism of $C^*L$-bimodules
  $C^*L\otimes (UL)^\vee\otimes C^*L= \left(C_*(L;UL)\otimes_{UL}C_*(UL;L)\right)^\vee \to C^*L$.
  The Hochschild cochain complexes of $UL$ and $C^*L$
  are defined by
   $$\hC^*(UL,UL) = \mbox{Hom}_{(UL)^e}(C_*(UL;L;UL), UL)$$
   $$\hC^*(C^*L,C^*L) =
   \mbox{Hom}_{(C^*L)^e}(\left(C_*(L;UL)\otimes_{UL}C_*(UL;L)\right)^\vee, C^*L)\,.$$
   The isomorphisms
   $$\mbox{Hom}_{(UL)^e}(C_*(UL;L;UL), UL) \cong Hom(C_*L,UL)$$
   $$\mbox{Hom}_{(C^*L)^e}(\left(C_*(L;UL)\otimes_{UL}C_*(UL;L)\right)^\vee,
   C^*L)\cong
   Hom((UL)^\vee, C^*L)$$
   induce  differentials on the right hand side terms.
   Now a simple computation (see for instance \cite{FMT}) shows that the map $D$
   that associates to a map in $\mbox{Hom}(C_*L,UL)$ the dual map in $\mbox{Hom}((UL)^\vee, C^*L)$
   is an isomorphism of complexes.
   This induces the following commutative diagram of complexes
   $$
   \begin{array}{ccc}
   UL & \stackrel{e'}{\longleftarrow} & \mbox{Hom}_{(UL)^e}(C_*(UL;L;UL), UL)\\
   \parallel && \parallel\\
   UL & \stackrel{e'}{\longleftarrow} & \mbox{Hom}(C_*L, UL)\\
   \downarrow {\scriptstyle\cong} && {\scriptstyle \cong}\downarrow {\scriptstyle D}\\
   \mbox{Hom}((UL)^\vee, \mathbb Q) &
   \stackrel{\mbox{\scriptsize Hom} ((UL)^\vee,\varepsilon)}{\longleftarrow}
   &\mbox{Hom}((UL)^\vee, C^*L)\\
   \parallel && \parallel\\
   \mbox{Hom}_{(C^*L)^e}(\left(C_*(L;UL)\otimes_{UL}C_*(UL;L)\right)^\vee, \mathbb Q)
   &\stackrel{\hC^*(C^*L,   \varepsilon)}
   {\longleftarrow} &
   \hC^*(C^*L,C^*L)\\
   \parallel &&\\
   \hC^*(C^*L,\mathbb Q) &
   &
   \end{array}
   $$
   where $\varepsilon : C^*L\to \mathbb Q$ denotes the canonical augmentation.
    The result follows now directly from the induced diagram in cohomology,
    combined with diagram (*).
   \hfill
   $\square$

 \vspace{1cm}
\hspace{-1cm}\begin{minipage}{19cm}
 \small
\begin{tabular}{lll}
felix@math.ucl.ac.be   &jean-claude.thomas@univ-angers.fr &
vigue@math.univ-paris13.fr\\
D\'epartement de math\'ematique  &  D\'epartement de
math\'ematique  &
D\'epartement de math\'ematique \\
 Universit\'e Catholique de Louvain  &Facult\'e des Sciences  &  Institut
Galil\'ee, \\
 2, Chemin du Cyclotron   &2, Boulevard Lavoisier &Universit\'e de
Paris-Nord\\
 1348 Louvain-La-Neuve, Belgium       & 49045 Angers, France & 93430
Villetaneuse, France
\end{tabular}

\end{minipage}

 \end{document}